\documentclass[journal]{IEEEtran}

\usepackage{graphicx}
\usepackage{ifpdf}
\usepackage{epstopdf}
\usepackage{ifpdf}
\ifpdf
\DeclareGraphicsExtensions{.eps}
\else
\DeclareGraphicsExtensions{.eps}
\fi
\usepackage{circuitikz}
\usetikzlibrary{arrows,matrix,backgrounds,calc,shapes,snakes}

\usepackage[noadjust]{cite}

\ifCLASSOPTIONcompsoc
\usepackage[caption=false,font=normalsize,labelfon
t=sf,textfont=sf]{subfig}
\else
\usepackage[caption=false,font=footnotesize]{subfi
	g}
\fi

\usepackage{eufrak}

\usepackage[cmex10]{amsmath}
\usepackage{amssymb}  
\usepackage{amsthm}
\usepackage{multirow}
\newtheorem{proposition}{Proposition}
\newtheorem{lemma}{Lemma}
\newtheorem{assumption}{Assumption}
\newtheorem{remark}{Remark}
\newtheorem{definition}{Definition}
\DeclareMathOperator*{\minimize}{minimize}

\DeclareMathOperator*{\subjectto}{subject~to}
\DeclareMathOperator*{\trace}{\mr{tr}}

\usepackage{steinmetz} 

\usepackage{threeparttable}
 \newcommand\MYhyperrefoptions{bookmarks=true,bookmarksnumbered=true,
	pdfpagemode={UseOutlines},plainpages=false,pdfpagelabels=true,
	colorlinks=true,linkcolor={black},citecolor={black},urlcolor={black},
	pdftitle={ThreePhaseVRs},
	pdfsubject={},
	pdfauthor={Mohammadhafez Bazrafshan, Nikolaos Gatsis, and Hao Zhu},
	pdfkeywords={}}

 \ifCLASSINFOpdf
\usepackage[\MYhyperrefoptions,pdftex]{hyperref}
\else
\usepackage[\MYhyperrefoptions,breaklinks=true,dvips]{hyperref}
\usepackage{breakurl}
\fi

\newcommand{\mr}[1]{\mathrm{#1}}
\newcommand{\mc}[1]{\mathcal{#1}}
\newcommand{\mb}[1]{\mathbf{#1}} 
\newcommand{\mbb}[1]{\mathbb{#1}}
\newcommand{\bmat}[1]{\begin{bmatrix} #1 \end{bmatrix}}
\newcommand{\smat}[1]{\left[\begin{smallmatrix} #1 \end{smallmatrix}\right]}
\newcommand{\diag}{\mathrm{diag}}
\newcommand{\mf}[1]{\mathfrak{#1}}
\renewcommand{\Re}[1]{\mr{Re}\left[{#1}\right]}
\renewcommand{\Im}[1]{\mr{Im}\left[{#1}\right]}

 \usepackage{algpseudocode}
\usepackage{algorithm}

\newcommand{\eq}[2]{
	\begin{IEEEeqnarray}{#1}
		#2
\end{IEEEeqnarray}}

	\definecolor{DarkBlue}{rgb}{0,0.1,0.7}

\begin{document}
	
	\title{Optimal Power Flow with Step-Voltage Regulators in Multi-Phase  Distribution Networks}
	\author{Mohammadhafez~Bazrafshan,~\IEEEmembership{ Member,~IEEE},~Nikolaos~Gatsis,~\IEEEmembership{Member,~IEEE},~and Hao~Zhu,~\IEEEmembership{Member,~IEEE}
		\thanks{ This material is based upon work supported by the National Science Foundation under Grant No. CCF-1421583, ECCS-1610732, and ECCS-1653706.}}

	\maketitle

	\begin{abstract}
	This paper develops a branch-flow based optimal power flow (OPF)  problem for multi-phase distribution networks that allows for tap selection of wye, closed-delta, and open-delta  step-voltage regulators (SVRs). SVRs are assumed ideal and their taps are represented by  continuous decision variables. To tackle the non-linearity, the branch-flow semidefinite programming framework of traditional OPF is expanded to accommodate SVR edges. Three types of non-convexity are addressed: (a) rank-1 constraints on non-SVR edges, (b) nonlinear equality constraints on SVR power flows and taps, and (c) trilinear equalities on SVR voltages and taps. 
	Leveraging a practical phase-separation assumption on the SVR secondary voltage, novel McCormick relaxations are provided for (c) and certain rank-1 constraints of (a), while dropping the rest. A linear relaxation based on conservation of power is used in place of (b).
  Numerical simulations on standard distribution test feeders corroborate the merits of the proposed convex formulation. 
	\end{abstract}


\section{Introduction}
\IEEEPARstart{T}{he} step-voltage regulator (SVR) is an autotransformer augmented by a tap-changing mechanism. It is used in medium-voltage distribution networks to maintain steady-state voltages within acceptable bounds.  Traditionally, SVR taps are automatically controlled via the line-drop compensator based on an approximate voltage-drop model from a local load-center~\cite{KerstingBook2001}. Such a scheme  is satisfactory for conventional distribution networks in which branch power flows are  unidirectional from the substation to the ends of the feeder.

Traditional  tap-selection is increasingly challenged by  modern distribution grids with high levels of distributed renewable generation.  A recent report~\cite{Nagarajan2018} highlights that during reversal of power flows,  the effectiveness of the regulator operation--as measured by voltage control per tap--reduces.  Of concern is also wear and tear of SVRs from  excessive tap changes following the fluctuation of renewables~\cite{Antoniadou-Plytaria2017,Agalgaonkar2014}.

Utilities with bundled retail and operations sectors can avoid  the aforementioned issues  by incorporating tap-selection into their optimal power flow (OPF) programs~\cite{Agalgaonkar2015}. Increasing renewable hosting capacity  by coordinating tap-selection and other voltage control resources is an additional advantage~\cite{Keane2011}.  Industry-provided integrated volt-var control applications for energy efficiency already support tap-setting of SVRs alongside with power factor and capacitor banks optimization; see e.g.,~the heuristic-based software product in~\cite{Yukon}.

Introducing SVR taps as variables  in distribution OPF is, however, technically challenging.  Primarily, they add to the non-convexity of the power flow equations. Since distribution networks are inherently unbalanced, tractable methodologies for multi-phase OPF problems  \cite{DallAnese2013,Gan2014,Gan2014Chordal, Zamzam2016, Wang2017a} should be expanded to this end.  Moreover, operational characteristics of various types of SVRs, i.e., wye, closed- or open-delta are different from each other. Since precise setting of SVRs aids in higher-quality voltage control,  raises the permissible loading level on feeders, and defers capacity investment costs~\cite{SystemsDivision}, a unified OPF program handling  various SVRs is  needed.

This paper develops an OPF that accounts for the tap selection of wye, closed-delta, and open-delta SVRs in multi-phase distribution networks.  To tackle the non-linearity, the branch-flow semidefinite programming (SDP) framework of multi-phase OPF is gracefully expanded to incorporate the full range of SVR models. Trilinear equalities that constrain SVR voltages and taps are handled via McCormick relaxations. The relaxation is made possible due to a phase-separation assumption on the SVR secondary voltage that is valid in practical multi-phase systems. This assumption is further leveraged  to approximate rank-1 constraints on the SVR secondary voltage matrix variable and improve the quality of the relaxation. The relevant literature is reviewed next.

\subsection{Literature review}
For single-phase radial networks, an OPF  considering tap-selection of on-load tap-changer transformers  is presented  in \cite{Wu2017} where the trilinear scalar constraint in transformer taps and voltages is converted to an exact mixed-binary linear constraint using binary expansion and big-M methods.  The second-order cone relaxation of  branch-flow power flows are then utilized to render an efficient mixed-integer second-order cone program (MISOCP).  An extension is presented in \cite{Chang2017} by incorporating static and discrete reactive power compensators.  

Considering unbalanced multi-phase operation, a  comprehensive OPF framework is developed in \cite{Paudyal2011}  by assembling ABCD parameters of transmission lines, transformers, as well as single- and three-phase wye tap-changers. Wye and delta-connected loads as well as switched shunt capacitors are also included.   The overall formulation is a mixed-integer nonlinear program (MINLP) which is then translated to a nonlinear program via a quadratic penalty function. Under a similar modeling framework, \cite{Daratha2014} also develops MINLP formulation of OPF to coordinate tap-changers and static var compensators with distributed generation which is ultimately solved via an ad-hoc two-stage procedure based on interior point branch and bound methods.   Despite their broad scope, these MINLP formulations turn out to be computationally intensive and may even yield locally suboptimal results. 

A more recent line of work explores convex relaxations. The work in  \cite{Robbins2016}  introduces the  tap selection of wye-connected SVRs inside the full SDP relaxation of the admittance-based OPF.    Power transfer from the primary to secondary of the SVR is accommodated by enforcing equal power injections on each side.  The trilinear matrix constraint in taps and voltages is further relaxed to a linear constraint that bounds the diagonals of the SDP variable using  minimum and maximum tap changes per phase.  By using the radial topology of distribution networks to improve computation time, \cite{Liu2017} leverages the chordal SDP relaxation of the admittance-based OPF. Further, the trilinear matrix constraint in taps and voltages is relaxed into  a linear semidefinite  matrix constraint that implicitly assumes that taps on every phase of the SVR are equal (gang-operated).  An unbalanced distribution reconfiguration problem is recently presented in \cite{Liu2018}, in  which SVR taps are represented via their binary expansion rendering a mixed-binary semidefinite program, albeit at the expense of introducing significant computational burden.  

The approaches in~\cite{Robbins2016,Liu2017, Liu2018} consider  wye-connected SVRs  for which the primary and secondary power injections are equal per  phase and the secondary voltage of each phase can be regulated independently from other phases.  Our previous work \cite{BazrafshanGatsisZhu2018} extends the chordal SDP relaxation of the admittance-based OPF to handle individually operated closed-delta and open-delta SVRs. 
However,  the formulation of~\cite{BazrafshanGatsisZhu2018} is only applicable to small-sized networks.

\subsection{Paper contributions and outline}
A  convex OPF formulation that can handle various types of SVRs and is applicable to larger networks is  missing in the literature. In this paper,  
	instead of using the admittance-based OPF \cite{Robbins2016,Liu2017,Liu2018,BazrafshanGatsisZhu2018}, the \emph{branch-flow} form of the power flow equations are leveraged to improve numerical stability \cite{Gan2014}.  Specifically, this work features the following contributions:
	\begin{itemize}
		\item A  branch-flow based OPF (BOPF)  is introduced that accommodates optimal tap selection of SVRs in multi-phase distribution networks. The formulation handles any combination of wye, closed-delta, and open-delta SVRs, as well as individual and gang operation of SVR taps.
		\item A nonlinear and non-convex SDP is developed that is provably equivalent to BOPF for radial networks, and a relaxation of BOPF for general meshed networks. The formulation extends the traditional framework of branch-flow SDP  for OPF originally put forth in~\cite{Gan2014} to incorporate most common SVR types.
		\item A novel convex relaxation is developed that ultimately alleviates the following non-convexity issues: (a) rank-1 constraints of non-SVR edges, (b) nonlinear constraints in SVR power flows and taps, and (c) trilinear constraints in SVR voltages and taps.  
	\end{itemize}

	The particular convex relaxation techniques are described next. Specifically, all rank-1 constraints of (a) are dropped except the ones that pertain to the SVR secondary, which are replaced by McCormick polyhedra. McCormick relaxations are also employed for (c), and a linear relaxation based on conservation of power is developed for (b). The resulting convex program is a tight relaxation of the original problem. 
The McCormick relaxations are enabled by a realistic assumption that the voltage angles on different phases of the SVR secondary are sufficiently separated. Different than this paper, McCormick relaxations for rank-1 constraints are adopted for single-phase networks~\cite{Kocuk2017}, by assuming phase differences between neighboring buses.

The proposed formulation is extensively tested on four standard distribution feeders that are properly edited to include wye, closed-delta, open-delta, and a mixture of SVRs. Detailed numerical comparisons with previously proposed convex techniques as well as with traditional nonlinear programming (NLP) algorithms are also provided. The findings indicate that the proposed convex formulation is capable of delivering tap settings of SVRs at almost zero optimality gaps (less than 1\%) in a time-span appropriate for OPF applications.

 The paper is organized as follows.  Notation,  network modeling including SVRs, and the non-convex OPF with SVRs are detailed in Section~\ref{Sec:NetModel}.  A rank-1 constrained OPF with SVRs   is introduced in Section~\ref{Sec:AltOPF} where the SVR non-convexities  represent themselves as trilinear equalities.  Convexifications of the SVR constraints as well as the rank-1 constraints via McCormick relaxations are pursued  in Section~\ref{Sec:ConvexOPF}. Formulation differences with prior work are highlighted in Section~\ref{Sec:LiteratureFormulations}.   Numerical tests that corroborate the practicality of the proposed formulation are carried out in Section~\ref{Sec:NumTests}.  The paper concludes in Section~\ref{Sec:Conclusions}.

 \section{Network Modeling and Branch-Flow OPF}
\label{Sec:NetModel}
This section introduces the notations and mathematical models for elements of the multi-phase distribution network including transmission lines, SVRs, and shunt elements.   The notation $\overline{(.)}$ is used to denote the conjugate transpose of $(.)$.

\subsection{General multi-phase notation}
A   multi-phase distribution network is mathematically modeled by a graph $(\mc{N}, \mc{E})$ where $\mc{N}$ is the set of buses and $\mc{E} \subseteq \mc{N} \times \mc{N} $ is the set of edges. The term ``edge'' is used instead of ``line'' to avoid  confusion. The set of buses represents shunt elements  and can be partitioned as $\mc{N}=\{0\} \cup \mc{N}_+$ where bus $0$ stands for the substation and the set $\mc{N}_+:=\{1,\ldots,N\}$ collects $N$ user buses. 

The set of edges $\mc{E}$  represents series elements of a distribution network and is partitioned as $\mc{E}=\mc{E}_{\mr{t}} \cup \mc{E}_{\mr{r}}$, where $\mc{E}_{\mr{t}}$ collects transmission lines and transformers while $\mc{E}_{\mr{r}}$ includes the SVRs.  An ordered pair $(n,m)$ (interchangeably, $n \rightarrow m$) belongs to the set $\mc{E}_{\mr{t}}$ when $n<m$ and bus $n$ is connected to bus $m$ via a transmission line or a transformer. An ordered pair $(n,m)$ belongs to the set $\mc{E}_{\mr{r}}$ when bus $n$ and $m$ are respectively  the primary of and secondary of an SVR.  The notation $n:n\rightarrow m$ means node $n \in \mc{N}$ such that $(n,m) \in \mc{E}$.  Define the set of primary nodes of SVRs connected to node $m$ as $\mc{N}_m^{\mr{p}}:=\{n: n \rightarrow m \in \mc{E}_{\mr{r}}\}$.

 The approach presented in this paper, as we will show in our numerical tests, is applicable to multi-phase networks with missing phases.   For the sake of exposition, however,  notations for strictly three-phase network are provided here. That is,  all buses and edges assume the phase set $\Omega= \{a,b,c\}$.  For $\phi \in \Omega$, denote the right shift as $\acute{a}=b$, $\acute{b}=c$, $\acute{c}=a$ and the left shift as $\grave{a}=c$, $\grave{b}=a$, $\grave{c}=b$. 

\subsection{Modeling of series elements}
\subsubsection{Transmission lines and transformers}
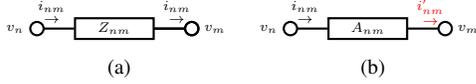
\begin{figure}[t]\scriptsize
	\centering
	\subfloat[]{\begin{tikzpicture}[scale=0.75] \tikzstyle{every node}=[scale=0.75]
		\tikzstyle{LineOut}   = [ draw=black,  thick, o-]
		\tikzstyle{LineIn}   = [ draw=black,  thick, -o]
		\tikzstyle{LittleArrow}= [draw=black, thin, ->]
		\tikzstyle{Block}=[rectangle, draw=black, thick, text width=5em, text centered, sharp corners]
		\node [matrix, ampersand replacement=\&, column sep=6mm, row sep=1mm] (SeriesModels){ 
			\node  (Noden) {$v_n$}; \& \node [Block] (Znm) {$Z_{nm}$}; \& \node (Nodem){$v_m$}; \\
		};
		\begin{scope}
		\path [LineOut](Noden.east) -- node [name=MidwayNode, midway, yshift=4mm] {$i_{nm}$}  (Znm.west); 
		\path [LineIn](Znm.east) -- node [name=MidwayNodeM, midway, yshift=4mm] {$i_{nm}$}  (Nodem.west); 
		\path [LineIn] (Znm.east) -- (Nodem.west);
		\draw ($(MidwayNode.south)$) node [yshift=-0.5mm] {$\rightarrow$}; 
		\draw ($(MidwayNodeM.south)$) node [yshift=-0.5mm] {$\rightarrow$}; 
		\end{scope}
		\end{tikzpicture} \label{Fig:TrLineDiagram}}{} 
	\subfloat[]{\begin{tikzpicture}[scale=0.75] \tikzstyle{every node}=[scale=0.75]
		\tikzstyle{LineOut}   = [ draw=black,  thick, o-]
		\tikzstyle{LineIn}   = [ draw=black,  thick, -o]
		\tikzstyle{LittleArrow}= [draw=black, thin, ->]
		\tikzstyle{Block}=[rectangle, draw=black, thick, text width=5em, text centered, sharp corners]
		\node [matrix, ampersand replacement=\&, column sep=6mm, row sep=1mm] (SeriesModels){ 
			\node  (Noden) {$v_n$}; \& \node [Block] (Anm) {$A_{nm}$}; \& \node (Nodem){$v_m$}; \\
		};
		\begin{scope}
		\path [LineOut](Noden.east) -- node [name=MidwayNodeN, midway, yshift=4mm] {$i_{nm}$}  (Anm.west); 
		\path [LineIn](Anm.east) -- node [name=MidwayNodeM, midway, yshift=4mm, color=red] {$i'_{nm}$}  (Nodem.west); 
		\draw ($(MidwayNodeN.south)$) node [yshift=-0.5mm] {$\rightarrow$}; 
		\draw ($(MidwayNodeM.south)$) node [color=red,yshift=-0.5mm] {$\rightarrow$}; 
		\end{scope}
		\end{tikzpicture} \label{Fig:SVRLineDiagram}}{}
	\caption{Series elements: \protect\subref{Fig:TrLineDiagram} Transmission lines and transformers
		\protect\subref{Fig:SVRLineDiagram} SVRs. For transformers and transmission lines on edge $(n,m)$, the per unit current flow from node $n$ to node $m$ is equal  as the per unit current flow \emph{received} at node $m$ from node $n$.  Therefore, only one variable $i_{nm}$ is required.  For SVRs,  the per unit current flow from node $n$ to node $m$ is related through~\eqref{Eq:CurrentGain} to the current flow \emph{received} at node $m$ from node $n$.  }
\end{figure}

Denote by $v_n, i_{nm} \in \mbb{C}^3$  and $Z_{nm} \in \mbb{C}^{3\times3}$ respectively the voltage phasor at node $n$, the  current phasor and the series impedance of the edge $(n,m) \in \mc{E}_t$ (see Fig.~\ref{Fig:TrLineDiagram}).  For wye-g--wye-g transformers, the series impedance is inverse of the per unit shunt admittance.  For other transformers, a suitable programming model would be to separate an invertible admittance from the common admittance models and reconnecting the remaining admittances as shunt (see e.g., \cite{Chen1991}).  Ohm's law implies
\eq{rCl}{
v_m=v_n-Z_{nm}i_{nm},  \: (n,m) \in \mc{E}_t. \label{Eq:Ohms}
}

\subsubsection{SVR modeling}
A three-phase SVR consists of three single-phase autotransformers that  typically connect in wye, closed-delta, or open-delta configuration.   The following modeling assumption on SVRs is asserted first. 
\begin{assumption}[Ideal SVRs] \label{Assumption:SVRs}
SVRs are ideal, i.e., the series impedance of the constituent autotransformers are negligible. 
\end{assumption}

Assumption~\ref{Assumption:SVRs} is realistic.  For instance,  \cite[Ch.~7]{KerstingBook2001} demonstrates that the per-unit series impedance of the autotransformer is approximately one tenth of that of the two-winding transformer and can be neglected for system-level studies. 

For edge $(n,m) \in \mc{E}_{\mr{r}}$, let $i_{nm}$ and $i'_{nm}$ respectively denote the  current phasors at primary and secondary of the SVR (see Fig.~\ref{Fig:SVRLineDiagram}). 
 Based on Assumption~\ref{Assumption:SVRs}, it suffices to model SVRs via their voltage and current gains as follows~\cite{BazrafshanGatsis2017}:
\eq{rCl}{
	\label{EqGroup:SVRGains}\IEEEyesnumber  \IEEEyessubnumber*
	v_n&=&A_{nm}v_m, \: (n,m) \in \mc{E}_{\mr{r}},  \label{Eq:VoltageGain} \IEEEyessubnumber\\
	i_{nm}&=&\bar{A}_{nm}^{-1}i'_{nm}, \:  (n,m) \in \mc{E}_{\mr{r}},\label{Eq:CurrentGain} \IEEEyessubnumber
}
where $A_{nm}$ is the voltage gain matrix and depends on the effective SVR turns ratio $r_{nm}$:
\eq{rCl}{
	A_{nm}= \diag(r_{nm}) D_{nm}+ F_{nm}, \: (n,m) \in \mc{E}_{\mr{r}}\label{Eq:ADecompose}
}
where  $D_{nm}$ and $F_{nm}$ are constant matrices given in Table~\ref{Table:SVRGains} for each SVR type. 
 For $(n,m) \in \mc{E}_{\mr{r}}$, the vector of effective turns ratios for wye, closed-delta, and open-delta SVRs is denoted by $r_{nm}:=\{r_{nm}^a, r_{nm}^b, r_{nm}^c\}$. For closed-delta SVRs, effective ratios on phase $ab$, $bc$, and $ca$ are given the labels $a$, $b$, and $c$. For open-delta SVRs,  effective ratio on phase $ab$ is given the label $a$ and effective ratio on phase $cb$ is given the label $c$. Open-delta SVRs do not have a third autotransformer, thus $r_{nm}^b=1$ is fixed and is not a variable.
 
  The relationship between the effective turns ratio and the taps for the SVR is 
 \eq{c}{
 	\mr{tap}^{\phi}=\mr{round}\left[\frac{\mp(1 -  r_{nm}^{\phi})}{0.00625}\right]. \label{Eq:Taps}
 }
 The  plus sign is used for type-A SVRs while the minus sign is used for type-B SVRs \cite{KerstingBook2001}.  The following modeling assumption regarding the SVR effective turns ratios is used for optimization.
 
 \begin{assumption}[Continuous turns ratios] \label{Assumption:SVRtaps}
	Effective turns ratios of SVRs  assume continuous values constrained by
	\eq{rCl}{
		r_{\min} \le r_{nm} \le r_{\max}, \: (n,m) \in \mc{E}_{\mr{r}} \label{Eq:RatioBounds}
	}
	where $[r_{\min},r_{\max}]=[0.9,1.1]$. 
\end{assumption}
This assumption is typical of works considering tap optimization of SVRs, see e.g.,~\cite{Robbins2016, Liu2017}. 
Under Assumption~\ref{Assumption:SVRtaps} holds, the taps span the interval $[-16,+16]$.   For open-delta SVRs, we set $r_{\min}^b=r_{\max}^b=1$ since it holds that $r_{nm}^b=1$. Assumptions~\ref{Assumption:SVRs} and~\ref{Assumption:SVRtaps} hold throughout this paper.

\begin{table}[t]
	\tiny
	\centering
	\caption{SVR voltage gain}
	\label{Table:SVRGains}
	\begin{IEEEeqnarraybox}[\IEEEeqnarraystrutmode ]{v.t/v/t/v/t/v/t.v}
			\IEEEeqnarraydblrulerowcut \\
					\IEEEeqnarrayseprow[2pt]\\
	&	SVR && $A_{nm}$ && $D_{nm}$ && $F_{nm}$ &\\
			\IEEEeqnarrayseprow[2pt]\\
		\IEEEeqnarraydblrulerowcut \\
			\IEEEeqnarrayseprow[5pt]\\
	&Wye && $\smat{ r_{nm}^a & 0 & 0 \\ 0 & r_{nm}^b & 0 \\ 0 & 0 & r_{nm}^c}$  && $\smat{1 &  0 & 0 \\ 0 & 1 & 0 \\ 0 & 0 & 1}$ && $\mb{O}$ & \\
	\IEEEeqnarrayseprow[5pt]\\
	\IEEEeqnarrayrulerow\\
		\IEEEeqnarrayseprow[5pt]\\
	&	Cl.-delta && $\smat{ r_{nm}^{a}  & 1- r_{nm}^{a} & 0 \\  0 & r_{nm}^{b} & 1- r_{nm}^{b} \\ 1-r_{nm}^{c} & 0 & r_{nm}^{c}}$ && $\smat{1 & -1 & 0 \\ 0 & 1 & -1\\ -1 & 0 & 1}$  && $\smat{ 0 & 1 & 0 \\ 0 & 0 & 1\\ 1 & 0 & 0}$ \\
		\IEEEeqnarrayseprow[5pt]\\
			\IEEEeqnarrayrulerow\\
					\IEEEeqnarrayseprow[5pt]\\
	&	Op.-delta &&  $\smat{ r_{nm}^{a} & 1-r_{nm}^{a} & 0 \\ 0 & 1 &  0 \\ 0 & 1-r_{nm}^{c} & r_{nm}^{c}}$ && $\smat{1 & -1 & 0 \\0 & 1 & 0 \\ 0 & -1  & 1}$ && $ \smat{0 & 1 & 0 \\ 0 & 0 & 0 \\ 0 & 1 & 0}$\\
			\IEEEeqnarrayseprow[5pt]\\
		\IEEEeqnarrayrulerow
		\end{IEEEeqnarraybox}
\end{table}

\color{black}
\subsection{Power balance equations}
The net current injection $i_m$ can be a sum of currents from a variety of sources. Here, we assume that the sources are  constant-power elements with complex power $\mf{s}$ as well as any constant-admittances (including capacitor banks and the sum of line shunt admittances) with admittance $Y_{m}$ connected at node $m$, as follows: 
\eq{rCl}{
	i_m=\diag(v_m^*)^{-1}\mf{s}_m^*-Y_m v_m, \: m \in \mc{N}.  \label{Eq:NetCurrentInjectionBreakDown}
}
Multiplying \eqref{Eq:NetCurrentInjectionBreakDown} by $\diag(v_m^*)$ and taking conjugate yields
\eq{rCl}{
	s_m=\mf{s}_m-\diag( v_m \bar{v}_m \bar{Y}_m), \: m \in \mc{N}. \label{Eq:NetPowerInjectionBreakDownNonConvex}
}
In~\eqref{Eq:NetPowerInjectionBreakDownNonConvex}, $s_m$ denotes the net complex power injection at node $m$, while $\mf{s}_m$ denotes the portion of the net complex power that originates from constant-power sources at node $m$.
Invoking KCL at bus $m$ yields 
\eq{rCl}{
	i_m&=&\sum\limits_{k:m\rightarrow k} i_{mk}-\sum\limits_{\substack{n:n\rightarrow m\\n \in \mc{N}_m^{\mr{p}}}}i'_{nm}- \sum\limits_{\substack{n:n\rightarrow m\\ n \notin \mc{N}_m^{\mr{p}} }}i_{nm}, \: m \in \mc{N}. \label{Eq:KCLatn} \IEEEeqnarraynumspace
}
Multiplying \eqref{Eq:KCLatn} by $\diag(v_m^*)$, taking conjugate, and again utilizing Lemma~\ref{Lemma:DiagToHermitian} yields
\eq{rCl}{
	s_m&=&\sum \limits_{k:m \rightarrow k} \diag(v_m \bar{i}_{mk})- \: \sum\limits_{\substack{n:n\rightarrow m\\n \in \mc{N}_m^{\mr{p}}}}\diag(v_m\bar{i}'_{nm})\notag \\
	&&- \sum \limits_{\substack{n:n\rightarrow m\\ n \notin \mc{N}_m^{\mr{p}} }}\diag(v_m\bar{i}_{nm}), \: m \in \mc{N}. \IEEEeqnarraynumspace \label{Eq:PowerBalanceNonConvex}
}
\subsection{Branch-flow optimal power flow with SVRs}
Let $s=\{s_m\}_{m \in \mc{N}}$, $\mf{s}=\{\mf{s}_m\}_{m \in \mc{N}}$ $v=\{v_m\}_{m \in \mc{N}}$, $i =\{i_{nm}\}_{(n,m) \in \mc{E}}$, and $i'=\{i'_{nm}\}_{(n,m) \in \mc{E}_\mr{r}}$, $r=\{r_{nm}\}_{(n,m) \in \mc{E}_{\mr{r}}}$, $A=\{A_{nm}\}_{(n,m) \in \mc{E}_{\mr{r}}}$. The branch flow formulation of optimal power flow problem (BOPF) with SVRs is given below: 
\eq{l}{\text{BOPF: }\label{EqGroup:OPF}\IEEEyesnumber \IEEEyessubnumber*
\minimize_{\substack{s,\mf{s}, v,i \\ i',r,A}}   c(s_0,\mf{s}, v, i)\hfill \label{Eq:OPFObj}\\
\subjectto \:  \eqref{Eq:Ohms}, \eqref{EqGroup:SVRGains},  \eqref{Eq:ADecompose}, \eqref{Eq:RatioBounds}, \eqref{Eq:NetPowerInjectionBreakDownNonConvex}, \eqref{Eq:PowerBalanceNonConvex} \notag \hfill \\
  v_0= \mr{v}_0 \hfill \label{Eq:Slack} \\
 v_{\min} \le |v_n| \le v_{\max}, \: n \in \mc{N} \label{Eq:VoltageBounds} \hfill\\
 \mf{s} \in \mc{S} \label{Eq:ConstantPowerSet} \hfill \IEEEeqnarraynumspace
}
where $\mr{v}_0$ is the fixed slack-bus voltage and \eqref{Eq:VoltageBounds} are the voltage limits. Equation \eqref{Eq:ConstantPowerSet} considers an operational set for constant-power injection. Usually, $\mc{S}=\prod_{m \in \mc{N}} \mc{S}_m$ where for distributed generation $\mc{S}_m$ is a disk while for constant-power loads, $\mc{S}_m$ is a singleton.  The cost, $c(s_0,\mf{s}, v,i)$ can account for thermal losses,  power import, or cost of distributed generation. 
 
The  BOPF formulation~\eqref{EqGroup:OPF}  incorporates  models of wye, closed-delta, and open-delta SVRs in the branch flow form of power flow equations.  BOPF is non-convex  due to bilinear and quadratic dependencies of \eqref{EqGroup:SVRGains}, \eqref{Eq:NetPowerInjectionBreakDownNonConvex}, and \eqref{Eq:PowerBalanceNonConvex} as well as the non-convexity imposed by the left-hand side of \eqref{Eq:VoltageBounds}.  BOPF is transformed in the next section to a rank-1 constrained nonlinear semidefinite program, which makes it amenable for branch-flow SDP relaxation.
\section{Rank-Constrained SDP For Branch-Flow OPF}
\label{Sec:AltOPF}
Let us introduce the following auxiliary matrix variables:
\eq{rCl}{
\label{EqGroup:Aux} \IEEEyesnumber \IEEEyessubnumber*
V_m &=&v_m\bar{v}_m, m \in \mc{N}  \IEEEyessubnumber \label{Eq:AuxV} \\
  I_{nm}&=&i_{nm}\bar{i}_{nm} (n,m) \in \mc{E}  \label{Eq:AuxI} \IEEEyessubnumber \\
 S_{nm}&=&v_{n} \bar{i}_{nm}, (n,m) \in \mc{E} \label{Eq:AuxPower} \IEEEyessubnumber \\
S'_{nm}&=&v_{m}\bar{i}'_{nm}, (n,m) \in \mc{E}_{\mr{r}}.   \label{Eq:AuxPrime}
}
Then,~\eqref{Eq:Ohms}, \eqref{EqGroup:SVRGains}, \eqref{Eq:NetPowerInjectionBreakDownNonConvex}, and \eqref{Eq:PowerBalanceNonConvex}  translate to 
\eq{rCl}{
V_m&=&V_n+Z_{nm}I_{nm}\bar{Z}_{nm} \hfill \notag \\
&& \: -(S_{nm}\bar{Z}_{nm}+Z_{nm}\bar{S}_{nm}), (n,m)  \in \mc{E}_{\mr{t}}. \IEEEeqnarraynumspace \hfill \label{Eq:OhmsMatrix}\\
V_n&=& A_{nm}V_m\bar{A}_{nm}, (n,m) \in \mc{E}_{\mr{r}} \label{Eq:VoltageGainTrilinear}  \IEEEeqnarraynumspace \hfill \\
0&=&\diag(A_{nm} ^{-1}S_{nm}A_{nm})- \diag(S'_{nm}), \:  (n,m) \in \mc{E}_{\mr{r}}\IEEEeqnarraynumspace \label{Eq:SVRPowers} \hfill \\
	s_m&=&\mf{s}_m-\diag(V_m\bar{Y}_m ), \: m \in \mc{N}. \label{Eq:NetPowerInjectionBreakDownConvex} \IEEEeqnarraynumspace \hfill \\
		s_m&=&\sum \limits_{k:m \rightarrow k} \diag(S_{mk})
	- \: \sum\limits_{\substack{n:n\rightarrow m\\n \in \mc{N}_m^{\mr{p}}}}\diag(S'_{nm}) \IEEEeqnarraynumspace \hfill \notag \\
	&& 	-\sum\limits_{\substack{n:n\rightarrow m\\ n \notin \mc{N}_m^{\mr{p}} }}\diag\left(S_{nm}-Z_{nm}I_{nm}\right), \: m \in \mc{N}.\label{Eq:PowerBalanceConvex}  \IEEEeqnarraynumspace \hfill 
}
where~\eqref{Eq:OhmsMatrix} and ~\eqref{Eq:VoltageGainTrilinear} are obtained by multiplying~\eqref{Eq:Ohms} and~\eqref{Eq:VoltageGain} by their  Hermitian. 
Equation~\eqref{Eq:SVRPowers} is obtained by multiplying~\eqref{Eq:VoltageGain} and Hermitian of~\eqref{Eq:CurrentGain},  incorporating~\eqref{Eq:AuxPower} and~\eqref{Eq:AuxPrime}, multiplying left and right respectively by $A_{nm}^{-1}$ and $A_{nm}$ and then taking only the diagonal elements.
Using~\eqref{Eq:AuxV} in~\eqref{Eq:NetPowerInjectionBreakDownNonConvex} yields~\eqref{Eq:NetPowerInjectionBreakDownConvex}.  Finally, using~\eqref{Eq:Ohms} to replace $v_m$ in the second line of~\eqref{Eq:PowerBalanceNonConvex} and subsequently substituting in~\eqref{EqGroup:Aux} yield~\eqref{Eq:PowerBalanceConvex}. Consider the following optimization problem:
\eq{l}{\text{RBOPF: } \label{EqGroup:OPFRank}\IEEEyesnumber \IEEEyessubnumber*
\minimize_{\substack{s,\mf{s}, V,I, S \\ S',r,A}}  c(s_0,\mf{s}, V,I) \IEEEeqnarraynumspace \hfill\label{Eq:OPFNonConvexObj}\\
	\subjectto \:   \eqref{Eq:ADecompose}, \eqref{Eq:RatioBounds}, \eqref{Eq:ConstantPowerSet}, \eqref{Eq:OhmsMatrix}, \eqref{Eq:VoltageGainTrilinear}, \eqref{Eq:SVRPowers}, \eqref{Eq:NetPowerInjectionBreakDownConvex}, \eqref{Eq:PowerBalanceConvex},  \notag \hfill \IEEEeqnarraynumspace\\
  V_0 = \mr{v}_0 \bar{\mr{v}}_0 \label{Eq:VSlack} \hfill \IEEEeqnarraynumspace \\	
 (v_{\min})^2 \le \diag(V_n) \le (v_{\max})^2, \: n \in \mc{N}  \label{Eq:VBounds} \hfill \IEEEeqnarraynumspace\\
\bmat{V_n & S_{nm} \\ \bar{S}_{nm} & I_{nm}} \succeq \mb{O}, \: (n,m) \in \mc{E} \IEEEeqnarraynumspace \label{Eq:PDConstraint} \hfill\\
\mr{rank}\left(\bmat{V_n & S_{nm} \\ \bar{S}_{nm} & I_{nm}}\right) =1, \: (n,m) \in \mc{E}.\IEEEeqnarraynumspace \label{Eq:Rank1Constraint} \hfill
}

The next two propositions characterize the relationship between  RBOPF and BOPF.
\begin{proposition}\label{Proposition:RBOPFRelax}
RBOPF is a relaxation of BOPF.
\end{proposition}
\begin{IEEEproof}
	If a point $\left(s,\mf{s}, v,i, i',r,A\right)$ is feasible for BOPF~\eqref{EqGroup:OPF}, then the point $\left(s,\mf{s}, V,I, S, S',r,A\right)$ obtained via~\eqref{EqGroup:Aux} is feasible for~\eqref{EqGroup:OPFRank}, as constraints~\eqref{Eq:OhmsMatrix}--\eqref{Eq:PowerBalanceConvex} together with \eqref{Eq:VSlack}--\eqref{Eq:Rank1Constraint} are satisfied.  The latter implies that the feasible set of RBOPF includes that of BOPF. 
	\end{IEEEproof}
The next proposition asserts that if the three-phase network has a radial topology, then RBOPF \eqref{EqGroup:OPFRank} is equivalent to BOPF~\eqref{EqGroup:OPF} by providing a unique way to go back from  $\left(s,\mf{s}, V,I, S, S',r,A\right)$ to $\left(s,\mf{s}, v,i, i',r,A\right)$.
\begin{proposition}\label{Proposition:Retrieve}
	Suppose the graph $(\mc{N},\mc{E})$ is radial and the point $\left(s,\mf{s}, V,I, S, S',r,A\right)$ is feasible for  \eqref{EqGroup:OPFRank}. Then, the point  $\left(s,\mf{s}, v,i, i',r,A\right)$ , where $v$, $i$, and $i'$ are computed via  Algorithm~\ref{Algorithm:Retrieve}, is feasible for~\eqref{EqGroup:OPF}. 
\end{proposition}
\begin{IEEEproof}
	The proof is provided in  Appendix \ref{Appendix:RetrieveProof}.  It relies on Lemma~\ref{Lemma:PowerSubstituteCurrent}, which states that conforming currents $i_{nm}$ and $i'_{nm}$ can be retrieved from RBOPF \eqref{EqGroup:OPFRank}.
	\end{IEEEproof}

\begin{lemma}\label{Lemma:PowerSubstituteCurrent}
	Suppose for $(v_n, v_m, i_{nm}, i'_{nm}, S_{nm}, S'_{nm})$ with $|v_m|>0$, equalities  \eqref{Eq:VoltageGain}, \eqref{Eq:AuxPower},  \eqref{Eq:SVRPowers} and the following hold:
		\eq{rCl}{
		\diag(S'_{nm})&=& \diag(v_{m}\bar{i}'_{nm}), \: (n,m) \in \mc{E}_{\mr{r}}. \label{Eq:DiagofAuxPower}
	}  Then, \eqref{Eq:CurrentGain} also holds.
\end{lemma}
\begin{IEEEproof}
	Substitute \eqref{Eq:AuxPower} and \eqref{Eq:DiagofAuxPower} into  \eqref{Eq:SVRPowers} to obtain
	\eq{rCl}{
		0&=& \diag(A_{nm} ^{-1}v_n\bar{i}_{nm}A_{nm})- \diag(v_m \bar{i}'_{nm}). \label{Eq:0DiagPowerSubCurrent}
	}
	Using  \eqref{Eq:VoltageGain} in \eqref{Eq:0DiagPowerSubCurrent} then yields
	\eq{rCl}{
		0&=&\diag\left[v_m(\bar{i}_{nm}A_{nm}-\bar{i}_{nm}')\right] \label{Eq:Proof1Diag1}. }
	Equation~\eqref{Eq:Proof1Diag1} is the pointwise multiplication of the non-zero vector $v_m$ with the vector $\bar{i}'_{nm}-\bar{i}_{nm}A_{nm}$.  Therefore, \eqref{Eq:CurrentGain} is inferred by  concluding that $\bar{i}'_{nm}-\bar{i}_{nm}A_{nm}=0$.
\end{IEEEproof}
 \begin{algorithm}[t] 
	\caption{Retrieve $v,i, i'$ from $V,I,S,S'$\label{Algorithm:Retrieve}}
\begin{algorithmic}[1]
		\State Initialize $\mc{N}_{\mr{visit}} := \{0\}$ and  $v_0=\mr{v}_0$.
		\While{$\mc{N}_{\mr{visit}}\neq \mc{N}$} 
		\State Find $(n,m) \in \mc{E}$ such that $n \in \mc{N}_{\mr{visit}}$ and $m \notin \mc{N}_{\mr{visit}}$. 
		\State	Set $i_{nm}:=\frac{1}{\trace{(V_n)}}\bar{S}_{nm}v_n \label{Eq:inmUpdate}$ 
		\If{$(n,m) \in \mc{E}_{\mr{r}}$}  
		\State Set $v_m:=A_{nm}^{-1} v_n \label{Eq:vmSVRUpdate}$ 
		\State Set $i'_{nm}:=\diag(v_m^*)^{-1} \diag(\bar{S}'_{nm}) \label{Eq:inmPrimeSVRUpdate}$
		\Else
		\State Set $v_m:=v_n- Z_{nm} i_{nm} \label{Eq:vmTrLineUpdate}$
		\EndIf
		\State Update $\mc{N}_{\mr{visit}}:=\mc{N}_{\mr{visit}} \cup \{m\}$.
		\EndWhile 
	\end{algorithmic}
\end{algorithm}
\begin{remark}
The radiality assumption in Proposition~\ref{Proposition:Retrieve} is leveraged only in the construction of Algorithm~\ref{Algorithm:Retrieve}, allowing for a way to compute a feasible point of BOPF from a feasible point of RBOPF. Equivalence between RBOPF and BOPF is thus only established for radial networks. However,  the ensuing convex relaxations for RBOPF are valid relaxations for BOPF under general network topologies,  as per Proposition~\ref{Proposition:RBOPFRelax}.
	\end{remark}
The  RBOPF~\eqref{EqGroup:OPFRank} is non-convex due to SVR constraints   \eqref{Eq:SVRPowers} and \eqref{Eq:VoltageGainTrilinear}   as well as the rank constraint~\eqref{Eq:Rank1Constraint}.   The next section examines convex alternatives for these constraints.

\section{Convex OPF with Tap Selection} 
\label{Sec:ConvexOPF}
\subsection{Convexifying the power equality \eqref{Eq:SVRPowers}}
\label{Subsec:ConvexPowerEquality}
Partition $\mc{E}_{\mr{r}}$  as $\mc{E}_{\mr{r}}=\mc{E}_{\mr{y}} \cup \mc{E}_{\mr{o}} \cup \mc{E}_{\mr{c}}$ where $\mc{E}_{\mr{y}}$, $\mc{E}_{\mr{o}}$, and $\mc{E}_{\mr{c}}$ respectively denote the set of wye, open-delta, and closed-delta SVRs. 
For wye SVRs, $A_{nm}$ is diagonal. Therefore, it is easily observed that constraint~\eqref{Eq:SVRPowers} is equivalent to
\eq{rCl}{
\diag(S_{nm})&=& \diag(S'_{nm}), \: (n,m) \in \mc{E}_{\mr{y}}. \label{Eq:WyeSVRPowers}
}
For open-delta and closed-delta SVRs, $A_{nm}$ is not diagonal and therefore \eqref{Eq:WyeSVRPowers} does not hold. In this case, due to the circular property of the trace of matrix products, we resort to the following \emph{relaxed} constraint on power conservation:
\eq{rCl}{
\trace{(S_{nm})}&=& \trace{(S'_{nm})}, \: (n,m) \in \mc{E}_{\mr{o}} \cup \mc{E}_{\mr{c}}. \label{Eq:DeltaSVRPowers}
}
\subsection{Convexifying the voltage equality \eqref{Eq:VoltageGainTrilinear}}\label{Subsec:ConvexVoltageEquality}
Define the following groups of variables:
\eq{rCl}{ \label{EqGroup:UWReImV} \IEEEyesnumber \IEEEyessubnumber*
	U_{n}&=&\mr{Re}[V_n], \:  W_{n} = \mr{Im}[V_n], \: n: (n,m) \in \mc{E}_{\mr{r}} \label{Eq:UnWn} \\
	U_{m}&=&\mr{Re}[V_m], \: W_{m} =\mr{Im}[V_m], \: m: (n,m) \in \mc{E}_{\mr{r}} \label{Eq:UmWm} \IEEEeqnarraynumspace \\
 \IEEEyesnumber \label{EqGroup:UmDUDDuF} \IEEEyessubnumber*
	\tilde{U}_{nm}&=&D_{nm} U_{m} \bar{D}_{nm}, \: (n,m) \in \mc{E}_{\mr{r}}  \label{Eq:UmDUD}\\
	\tilde{W}_{nm}&=& D_{nm}W_{m}\bar{D}_{nm},  \: (n,m) \in \mc{E}_{\mr{r}} \label{Eq:WmDWD} \\
		\hat{U}_{nm} &=& D_{nm}U_m\bar{F}_{nm}, \: (n,m) \in \mc{E}_{\mr{r}} \label{Eq:UmDUF} \\
	\hat{W}_{nm}&=&D_{nm}W_{m}\bar{F}_{nm}, \: (n,m) \in \mc{E}_{\mr{r}} \label{Eq:WmDWF} \\
 \IEEEyesnumber \label{EqGroup:UnRUR} \IEEEyessubnumber* \IEEEeqnarraynumspace
	\tilde{\mf{U}}_{nm} &=& \diag(r_{nm}) \tilde{U}_{nm} \diag(r_{nm}), \: (n,m) \in \mc{E}_{\mr{r}}   \label{Eq:UnRUR} \\
	\tilde{\Psi}_{nm} &=& \diag(r_{nm}) \tilde{W}_{nm} \diag(r_{nm}), \: (n,m) \in \mc{E}_{\mr{r}}   \label{Eq:WnRUR} \\
		\hat{\mf{U}}_{nm}&=&\diag(r_{nm}) \hat{U}_{nm}, \: (n,m) \in \mc{E}_{\mr{r}} \label{Eq:UnRU} \\
	\hat{\Psi}_{nm}&=&\diag(r_{nm}) \hat{W}_{nm}, \: (n,m) \in \mc{E}_{\mr{r}}. \label{Eq:WnRW} 
}

Using~\eqref{Eq:ADecompose} and~\eqref{EqGroup:UWReImV}--\eqref{EqGroup:UnRUR}, constraint~\eqref{Eq:VoltageGainTrilinear} is recast as
\eq{rCl}{ \label{EqGroup:LinearTrilinear} \IEEEyesnumber \IEEEyessubnumber*
U_n &=& \tilde{\mf{U}}_{nm}+\hat{\mf{U}}_{nm}+\bar{\hat{\mf{U}}}_{nm}+F_{nm}U_m \bar{F}_{nm}, (n,m) \in \mc{E}_{\mr{r}} \label{Eq:UnUmFUF} \IEEEeqnarraynumspace\\
W_n &=& \tilde{\Psi}_{nm}+\hat{\Psi}_{nm}-\bar{\hat{\Psi}}_{nm}+F_{nm}W_m \bar{F}_{nm}, (n,m) \in \mc{E}_{\mr{r}}. \label{Eq:WnWmFWF} \IEEEeqnarraynumspace
}
The nonconvexity now lies only in \eqref{EqGroup:UnRUR}.   Based on Hermitian symmetry of $V_n$ and $V_m$, and recalling that $\acute{\phi}$ is the right shift of phase $\phi$, 
\eqref{EqGroup:UnRUR} is equivalent to
 \eq{rCll}{ \label{EqGroup:UnWnr2UmWm} \IEEEyesnumber \IEEEyessubnumber*
 \tilde{\mf{U}}_{nm}^{\phi\phi'} &=& r_{nm}^{\phi}r_{nm}^{\phi'} \tilde{U}_{nm}^{\phi \phi'
 }, \: &\phi \in \Omega, \phi' \in \{\phi, \acute{\phi}\}, \IEEEeqnarraynumspace \label{Eq:Unr2Um}  \\
 \tilde{\Psi}_{nm}^{\phi\phi'} &=& r_{nm}^{\phi}r_{nm}^{\phi'} \tilde{W}_{nm}^{\phi \phi'
},  \: &\phi \in \Omega, \phi' \in \{\phi, \acute{\phi}\}, \IEEEeqnarraynumspace \label{Eq:Wnr2Wm} \\
 \hat{\mf{U}}_{nm}^{\phi\phi'} &=& r_{nm}^{\phi} \hat{U}_{nm}^{\phi \phi'
}, \: &\phi , \phi' \in \Omega,  \IEEEeqnarraynumspace \label{Eq:UnrUm} \\
 \hat{\Psi}_{nm}^{\phi\phi'} &=& r_{nm}^{\phi} \hat{W}_{nm}^{\phi \phi'
},  \: &\phi, \phi' \in \Omega. \IEEEeqnarraynumspace \label{Eq:WnrWm}}

Linear relaxations of the bilinear and trilinear equalities in~\eqref{EqGroup:UnWnr2UmWm} are based on McCormick envelopes given in Definition~\ref{Definition:McCormick}.   To employ McCormick envelopes,  bounds on $r$ and $v$ provided in \eqref{Eq:RatioBounds} and \eqref{Eq:VoltageBounds} are leveraged together with the following assumption on the secondary voltage of the SVR. 
\begin{assumption}[Phase separation]
	\label{Assumption:PhaseSeparation}
	Let the complex voltage phasor on the secondary of an SVR be equal to $v_m=\{|v_m^a|e^{j\theta^a}, |v_m^b| e^{j \theta^b}, |v_m^c| e^{j \theta^c}\}$ for $m: n \rightarrow m$. There exists $\Delta >0$ such that $\theta^a$, $\theta^b$, and $\theta^c$ satisfy
	\eq{rCl}{
		90^\circ \le 120- \Delta \le \theta^{\phi}-\theta^{\acute{\phi}} \le 120+ \Delta \le 180^\circ, \phi \in \Omega. \label{Eq:ThetaBounds} \IEEEeqnarraynumspace
	}
\end{assumption}
Assumption~\ref{Assumption:PhaseSeparation} is based on the fact that phases of a distribution network, are well separated even under unbalanced operation. Based on Assumption~\ref{Assumption:PhaseSeparation}, the following proposition is provided whose proof is furnished in Appendix~\ref{Appendix:ProofUWBounds}.
\begin{proposition}\label{Proposition:UWBounds}
	Under Assumption~\ref{Assumption:PhaseSeparation} and the bounds in~\eqref{Eq:VoltageBounds}, entry-wise bounds on  $U_m$, $W_m$, $\tilde{U}_{nm}$, $\tilde{W}_{nm}$, $\hat{U}_{nm}$, and $\hat{W}_{nm}$ for $m:(n,m) \in \mc{E}_{\mr{r}}$ are computed given $v_{\min}$, $v_{\max}$, and $\Delta$:
	\eq{rCl}{ \IEEEyesnumber  \label{EqGroup:UWBounds}
		U_{\min} &\le& U_m \le U_{\max}, \: m: (n,m) \in \mc{E}_{\mr{r}}  \IEEEeqnarraynumspace\label{Eq:UBounds} \IEEEyessubnumber \\
		W_{\min} &\le& W_m \le W_{\max},\:  m:(n,m) \in \mc{E}_{\mr{r}} \label{Eq:WBounds} \IEEEyessubnumber \\ 
	 \IEEEyesnumber  \label{EqGroup:TildeUWBounds}
	 	\tilde{U}_{\min} &\le& \tilde{U}_{nm} \le \tilde{U}_{\max}, \: m: (n,m) \in \mc{E}_{\mr{r}} \IEEEyessubnumber \label{Eq:TildeUBounds} \\ 
		\tilde{W}_{\min} &\le& \tilde{W}_{nm} \le\tilde{W}_{\max}, \: m: (n,m) \in \mc{E}_{\mr{r}}  \IEEEyessubnumber \label{Eq:TildeWBounds} \\
		 \IEEEyesnumber  \label{EqGroup:HatUWBounds}
		 	\hat{U}_{\min} &\le& \hat{U}_{nm} \le \hat{U}_{\max}, \: m: (n,m) \in \mc{E}_{\mr{r}} \IEEEyessubnumber \label{Eq:HatUBounds}  \\
		\hat{W}_{\min} &\le& \hat{W}_{nm} \le \hat{W}_{\max},\:  m: (n,m) \in \mc{E}_{\mr{r}}.  \IEEEyessubnumber \label{Eq:HatWBounds}  	
	}
\end{proposition}

\begin{definition} \label{Definition:McCormick}
For variables $u$, $w$, and $x$ as well as the given parameters $u_{\min}$, $u_{\max}$, $w_{\min}$, and $w_{\max}$ with $u_{\min} \le u_{\max}$ and $w_{\min}\le w_{\max}$, consider the following set of inequalities:
\eq{rCl}{
	\IEEEyesnumber \label{EqGroup:uwxMcCormick}  \IEEEyessubnumber* 
	u_{\min} \le u \le u_{\max} \label{Eq:uBounds} \IEEEeqnarraynumspace\\
	w_{\min} \le w \le w_{\max} \label{Eq:wBounds} \IEEEeqnarraynumspace \\
	u_{\min}w+uw_{\min}-u_{\min}w_{\min} &\le& x \IEEEeqnarraynumspace\label{Eq:uwxUnderEstimatorMin} \\
	u_{\max}w+uw_{\max}-u_{\max}w_{\max} &\le& x \IEEEeqnarraynumspace\label{Eq:uwxUnderEstimatorMax} \\
	u_{\max}w+uw_{\min}-u_{\max}w_{\min} &\ge& x \IEEEeqnarraynumspace\label{Eq:uwxOverEstimatorMaxMin} \\
	u_{\min}w+uw_{\max} - u_{\min}w_{\max} &\ge& x\IEEEeqnarraynumspace \label{Eq:uwxOverEstimatorMinMax}
}
We compactly denote \eqref{EqGroup:uwxMcCormick} by
\eq{rCl}{
\mc{M}\left(u,w,x; u_{\min},u_{\max}, w_{\min},w_{\max}\right) \le 0 \label{Eq:uwxM}
} We refer to \eqref{Eq:uwxM} as the McCormick polyhedron of variables $u$, $w$, and $x$, which is a linear relaxation of the bilinear constraint $x=uw$ when $u$ and $w$ are bounded by \eqref{Eq:uBounds} and \eqref{Eq:wBounds}.
\end{definition}
Let us introduce the additional variables $R_{nm}^{\phi \phi'}$ constrained as follows:
	\eq{rCl}{R_{nm}^{\phi \phi'}=r_{nm}^{\phi}r_{nm}^{\phi'}, (n,m) \in \mc{E}_{\mr{r}}, \phi,\phi' \in \Omega. \label{Eq:R}}
 The bounds in~\eqref{Eq:RatioBounds} enable the following relaxation for~\eqref{Eq:R}:
	\eq{rCl}{(n,m) \in \mc{E}_{\mr{r}}, \phi,\phi' \in \Omega: \hfill \notag \\
		\mc{M}\left(r_{nm}^{\phi}, r_{nm}^{\phi'}, R_{nm}^{\phi \phi'}; r_{\min}, r_{\max}, r_{\min}, r_{\max}\right) &\le& 0. \IEEEeqnarraynumspace \label{Eq:rMcCormick} 
	}%
Further, notice that \eqref{Eq:R} is equivalent to $R_{nm}=r_{nm}\bar{r}_{nm}$ where matrix $R_{nm}$ is assembled by concatenating the values of $R_{nm}^{\phi,\phi'}$ for $\phi,\phi' \in \Omega$. Therefore, a semidefinite relaxation of~\eqref{Eq:R} may be additionally used:
\eq{rCl}{ \bmat{R_{nm} & r_{nm} \\ \bar{r}_{nm} & 1} \succeq \mb{O}, (n,m) \in \mc{E}_{\mr{r}}.  \label{Eq:RPD}}
Lower and upper bounds on variables $\tilde{U}_{nm}$, $\hat{U}_{nm}$, $\tilde{W}_{nm}$, and $\hat{W}_{nm}$ are provided by Proposition~\ref{Proposition:UWBounds}. Upon substituting~\eqref{Eq:R} into~\eqref{Eq:Unr2Um} and~\eqref{Eq:Wnr2Wm} and utilizing the bounds of Proposition~\ref{Proposition:UWBounds}, the constraints in~\eqref{EqGroup:UnWnr2UmWm} are respectively relaxed to \eq{rCl}{ \IEEEyesnumber \label{EqGroup:UnWnr2UmWmMcCormick}  \IEEEyessubnumber* (n,m)\in \mc{E}_{\mr{r}}, \phi \in \Omega,  \phi' \in \{\phi, \acute{\phi}\}:\hfill  \notag \\
\mc{M}\left(R_{nm}^{\phi\phi'}, \tilde{U}_{nm}^{\phi\phi'}, \tilde{\mf{U}}_{nm}^{\phi \phi'}; r_{\min}^2, r_{\max}^2, \tilde{U}_{\min}^{\phi \phi'}, \tilde{U}_{\max}^{\phi \phi'}\right) &\le& 0\IEEEeqnarraynumspace \label{Eq:TildeUr2McCormick} \\
\mc{M}\left(R_{nm}^{\phi\phi'}, \tilde{W}_{nm}^{\phi\phi'}, \tilde{\Psi}_{nm}^{\phi \phi'}; r_{\min}^2, r_{\max}^2, \tilde{W}_{\min}^{\phi \phi'}, \tilde{W}_{\max}^{\phi \phi'}\right) &\le& 0\IEEEeqnarraynumspace \label{Eq:TildeWr2McCormick} \\
(n,m) \in \mc{E}_{\mr{r}}, \phi, \phi' \in \Omega: \hfill \notag \\
\mc{M}\left(r_{nm}^{\phi}, \hat{U}_{nm}^{\phi\phi'}, \hat{\mf{U}}_{nm}^{\phi \phi'}; r_{\min}, r_{\max}, \hat{U}_{\min}^{\phi \phi'}, \hat{U}_{\max}^{\phi \phi'}\right) &\le& 0\IEEEeqnarraynumspace \label{Eq:HatUrMcCormick}  \\
\mc{M}\left(r_{nm}^{\phi}, \hat{W}_{nm}^{\phi\phi'}, \hat{\Psi}_{nm}^{\phi \phi'}; r_{\min}, r_{\max}, \hat{W}_{\min}^{\phi \phi'}, \hat{W}_{\max}^{\phi \phi'}\right) &\le& 0.\IEEEeqnarraynumspace \label{Eq:HatWrMcCormick}}
\subsection{Rank reinforcements} \label{Subsec:RankReinforcements}
Recall that the third source of nonconvexity in  RBOPF~\eqref{EqGroup:OPFRank} is the Rank-1 constraint \eqref{Eq:Rank1Constraint}. The goal here is to improve the quality of the voltage solution provided by the relaxation of RBOPF by approximating the constraint \vspace{-0.1cm}
\eq{rCl}{
	\mr{Rank}(V_m)&=&1 \label{Eq:VnRank1} \vspace{-0.1cm}
}which is a consequence of~\eqref{Eq:AuxV}.  We first borrow the following result \cite[Proposition~3.1]{Kocuk2017}.  
\begin{proposition}\label{Proposition:MatrixMinors}
	The Hermitian matrix $V_m$ is positive semidefinite and rank-1 if and only if the diagonal entries of $V_m$ are nonnegative and all of  $2\times 2$ minors of $V_m$ are zero.
\end{proposition} \vspace{-0.1cm}
We use Propositions~\ref{Proposition:UWBounds} and Proposition~\ref{Proposition:MatrixMinors} to provide a linear relaxation of \eqref{Eq:VnRank1}.  Since $V_{m} \in \mbb{C}^{3 \times 3}$ is  Hermitian, setting its minors to zero yields $9$ equalities: 
\eq{rCl}{ \IEEEyesnumber \label{EqGroup:Minors} \IEEEyessubnumber*
U_m^{\phi\phi} U_m^{\acute{\phi}\acute{\phi}}-\left(U_m^{\phi\acute{\phi}}\right)^2-\left(W_{m}^{\phi\acute{\phi}}\right)^2&=&0, \phi \in  \Omega \label{Eq:PrincipalMinorForm} \\
U_{m}^{\grave{\phi} \phi} U_{m}^{\phi \acute{\phi}} - W_m^{\phi \acute{\phi}} W_{m}^{\grave{\phi} \phi} - U_m^{\phi \phi} U_{m}^{\grave{\phi} \acute{\phi}} &=&0, \phi \in \Omega  \label{Eq:MinorRealForm}\\
 -U_{m}^{\grave{\phi} \phi} W_m^{\phi \acute{\phi}}+U_{m}^{ \phi \acute{\phi}} W_{m}^{\grave{\phi} \phi} - U_m^{\phi \phi} W_{m}^{\grave{\phi} \acute{\phi}} &=&0, \phi \in \Omega.  \label{Eq:MinorImagForm}
\vspace{-0.1cm}} Define the following variables for  $m: (n,m) \in \mc{E}_{r},\phi \in \Omega$: 
\eq{rClrClrCl}{
\IEEEyesnumber \label{EqGroup:Xm1To9}  \IEEEyessubnumber* 
X_m^{\phi1}&=&U_m^{\phi\phi} U_m^{\acute{\phi}\acute{\phi}}, \:  & X_{m}^{\phi2}&=&\left(U_m^{\phi \acute{\phi}}\right)^2, \: & X_{m}^{\phi3}&=&\left(W_m^{\phi \acute{\phi}}\right)^2 \IEEEeqnarraynumspace \label{Eq:Xm1To3}\\
X_m^{\phi4}&=&U_m^{\grave{\phi}\phi} U_m^{\phi\acute{\phi}}, \:  & X_{m}^{\phi5}&=& W_m^{\phi \acute{\phi}}W_m^{\grave{\phi} \phi} , \: & X_{m}^{\phi6}&=&U_m^{\phi \phi} U_m^{\grave{\phi} \acute{\phi}} \IEEEeqnarraynumspace \label{Eq:Xm4To6} \\
X_{m}^{\phi7}&=&U_m^{\grave{\phi} \phi}  W_m^{\phi \acute{\phi}}, \:  &
X_m^{\phi7}&=&U_m^{\phi \acute{\phi}} W_m^{\grave{\phi}\phi},  \: & X_{m}^{\phi9}&=&U_m^{\phi \phi} W_m^{\grave{\phi} \acute{\phi}}. \IEEEeqnarraynumspace  \label{Eq:Xm7To9}
} 
By capturing the bilinear relation in \eqref{EqGroup:Minors} using McCormick envelopes, we can again obtain its linear relaxation as  
{\small \eq{rCl}{\IEEEyesnumber \label{EqGroup:MinorsLinear}  \IEEEyessubnumber*
X_{m}^{\phi1}- X_{m}^{\phi2}- X_m^{\phi3}&=&0\label{Eq:PrincipalMinorLinear}  \IEEEeqnarraynumspace\\
X_{m}^{\phi4}-X_{m}^{\phi5}-X_m^{\phi6}&=&0 \IEEEeqnarraynumspace\label{Eq:MinorRealLinear}\\
-X_{m}^{\phi7}+X_{m}^{\phi8}-X_m^{\phi9}&=&0 \IEEEeqnarraynumspace\label{Eq:MinorImagLinear} \vspace{-0.3cm} } 
 \eq{rCl}{  \IEEEyessubnumber* 
	\mc{M}\left(U_m^{\phi\phi}, U_m^{\acute{\phi}\acute{\phi}}, X_m^{\phi1}; U_{\min}^{\phi \phi}, U_{\max}^{\phi \phi}, U_{\min}^{\acute{\phi}\acute{\phi}}, U_{\max}^{\acute{\phi}\acute{\phi}}\right) &\le& 0\IEEEeqnarraynumspace \label{Eq:Xm1McCormick} \\
	\mc{M}\left(U_m^{\phi \acute{\phi}},U_m^{\phi \acute{\phi}}, X_m^{\phi2}; U_{\min}^{\phi \acute{\phi}}, U_{\max}^{\phi \acute{\phi}}, U_{\min}^{\phi \acute{\phi}}, U_{\max}^{\phi \acute{\phi}}\right) &\le& 0\IEEEeqnarraynumspace \label{Eq:Xm2McCormick} \\
	\mc{M}\left(W_m^{\phi \acute{\phi}},W_m^{\phi \acute{\phi}}, X_m^{\phi3}; W_{\min}^{\phi \acute{\phi}}, W_{\max}^{\phi \acute{\phi}}, W_{\min}^{\phi \acute{\phi}}, W_{\max}^{\phi \acute{\phi}}\right) &\le& 0\IEEEeqnarraynumspace \label{Eq:Xm3McCormick} \\
	\mc{M}\left(U_m^{\grave{\phi} \phi}, U_m^{\phi \acute{\phi}}, X_m^{\phi4}; U_{\min}^{\grave{\phi} \phi}, U_{\max}^{\grave{\phi} \phi}, U_{\min}^{\phi \acute{\phi}}, U_{\max}^{\phi \acute{\phi}}\right) &\le& 0\IEEEeqnarraynumspace \label{Eq:Xm4McCormick} \\
	\mc{M}\left(W_m^{\grave{\phi} \phi}, W_m^{\phi \acute{\phi}}, X_m^{\phi5}; W_{\min}^{\grave{\phi} \phi}, W_{\max}^{\grave{\phi} \phi}, W_{\min}^{\phi \acute{\phi}}, W_{\max}^{\phi \acute{\phi}}\right) &\le& 0\IEEEeqnarraynumspace \label{Eq:Xm5McCormick} 			\\
	\mc{M}\left(U_m^{\phi \phi}, U_m^{\grave{\phi} \acute{\phi}}, X_m^{\phi6}; U_{\min}^{\phi\phi}, U_{\max}^{\phi \phi}, U_{\min}^{\grave{\phi} \acute{\phi}}, U_{\max}^{\grave{\phi}\acute{\phi}}\right) &\le& 0\IEEEeqnarraynumspace \label{Eq:Xm6McCormick} \\	
	\mc{M}\left(U_m^{\grave{\phi} \phi}, W_m^{\phi \acute{\phi}}, X_m^{\phi7}; U_{\min}^{\grave{\phi} \phi}, U_{\max}^{\grave{\phi} \phi}, W_{\min}^{\phi \acute{\phi}}, W_{\max}^{\phi \acute{\phi}}\right) &\le& 0\IEEEeqnarraynumspace \label{Eq:Xm7McCormick} \\
		\mc{M}\left(U_m^{\grave{\phi} \phi}, W_m^{\phi \acute{\phi}}, X_m^{\phi8}; U_{\min}^{\grave{\phi} \phi}, U_{\max}^{\grave{\phi} \phi}, W_{\min}^{\phi \acute{\phi}}, W_{\max}^{\phi \acute{\phi}}\right) &\le& 0\IEEEeqnarraynumspace \label{Eq:Xm8McCormick} 			\\	
	\mc{M}\left(U_m^{\phi \phi}, W_m^{\grave{\phi} \acute{\phi}}, X_m^{\phi9}; U_{\min}^{\phi\phi}, U_{\max}^{\phi \phi}, W_{\min}^{\grave{\phi} \acute{\phi}}, W_{\max}^{\grave{\phi}\acute{\phi}}\right) &\le& 0\IEEEeqnarraynumspace \label{Eq:Xm9McCormick} 	
	}} for $m: (n,m) \in \mc{E}_{\mr{r}}, \phi \in \Omega$.
\subsection{Convex relaxation of BOPF with SVRs}
\label{Subsec:BOPF}
\begin{figure}[t]\scriptsize
	\centering
	\subfloat[]{\begin{tikzpicture}[scale=1]
\node[circle, draw, fill=gray!5 ,text height = 1.5cm,minimum     width=2cm] at (0,0) (MBOPF){};
\draw node at ($(MBOPF.south)+(0,0.2)$) {MBOPF};
\node[star, star points=7, star point ratio=0.85, draw, fill=gray!20 ,text height = .35cm,minimum     width=1.3cm] at  ($(MBOPF.north)-(0,0.8)$) (RBOPF){};
\draw node at ($(RBOPF.west)+(0.7,-0.4)$) {RBOPF};
\node[star, star points=9, star point ratio=0.85, draw, fill=gray!50 ,text height = .19cm,minimum     width=0.75cm] at  ($(RBOPF.north)-(0,0.5)$) (BOPF){};
\draw node at ($(BOPF.center)$) {BOPF};
		\end{tikzpicture} \label{Fig:FeasibleMeshed}}{}  	 \hspace{1.5cm}
	\subfloat[]{\begin{tikzpicture}[scale=1]
		\node[circle, draw, fill=gray!5 ,text height = 1.5cm,minimum     width=2cm] at (0,0) (MBOPF){};
		\draw node at ($(MBOPF.south)+(0,0.2)$) {MBOPF};
		\node[star, star points=9, star point ratio=0.85, draw, fill=gray!20 ,text height = .35cm,minimum     width=1.3cm] at  ($(MBOPF.north)-(0,0.8)$) (RBOPF){};
				\draw node at ($(RBOPF.north)+(0,-0.5)$) {BOPF};
			\draw node at ($(RBOPF.north)+(0,-0.75)$) {\&};
		\draw node at ($(RBOPF.north)+(0,-1)$) {RBOPF};
		\end{tikzpicture} \label{Fig:FeasibleRadial}}{}
	\caption{Schematic representation of the feasible sets for BOPF, RBOPF, and MBOPF: \protect\subref{Fig:FeasibleMeshed} Meshed networks;
		\protect\subref{Fig:FeasibleRadial} Radial networks. \label{Fig:FeasibleSets} }
\end{figure}
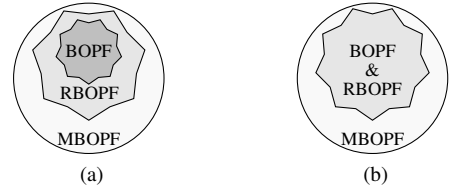

The proposed convex formulation, MBOPF, is 
\eq{l}{\text{MBOPF:} \label{EqGroup:MBOPF}\IEEEyesnumber 
	\minimize_{\substack{s,\mf{s}, V,I, S \\ S',r,R,A, U \\ W,  \tilde{U}, \tilde{W}, \hat{U}, \hat{W} \\ \tilde{\mf{U}}, \tilde{\Psi}, \hat{\mf{U}},  \hat{\Psi}, X}}  c(s_0,\mf{s}, V,I) \IEEEeqnarraynumspace \hfill\label{Eq:MBOPFObj}\\
	\subjectto   \eqref{Eq:ADecompose}, \eqref{Eq:RatioBounds},  \eqref{Eq:ConstantPowerSet}, \eqref{Eq:OhmsMatrix}, \eqref{Eq:NetPowerInjectionBreakDownConvex}, \eqref{Eq:PowerBalanceConvex},  \hfill \notag \\
 \eqref{Eq:VSlack}, \eqref{Eq:VBounds},  \eqref{Eq:PDConstraint}, \eqref{Eq:WyeSVRPowers}, \eqref{Eq:DeltaSVRPowers}, \hfill \notag \\
 \eqref{EqGroup:UWReImV}, \eqref{EqGroup:UmDUDDuF}, \eqref{EqGroup:LinearTrilinear} , 
\eqref{EqGroup:UWBounds}, \eqref{EqGroup:TildeUWBounds}, \eqref{EqGroup:HatUWBounds}, \hfill \notag \\ \eqref{Eq:rMcCormick}, \eqref{Eq:RPD}, \eqref{EqGroup:UnWnr2UmWmMcCormick}, \eqref{EqGroup:MinorsLinear}.  \notag \hfill  \IEEEeqnarraynumspace 
}
The following proposition clarifies the relationship between MBOPF and RBOPF.
\begin{proposition} \label{Proposition:MBOPF}
	Under Assumption~\ref{Assumption:PhaseSeparation}, MBOPF is a relaxation of RBOPF. 
\end{proposition}
\begin{IEEEproof}
Constraint~\eqref{Eq:SVRPowers} of RBOPF is relaxed to constraints~\eqref{Eq:WyeSVRPowers} and~\eqref{Eq:DeltaSVRPowers} of MBOPF.  Constraint~\eqref{Eq:VoltageGainTrilinear} of RBOPF is relaxed to constraints  \eqref{EqGroup:UWReImV},  \eqref{EqGroup:UmDUDDuF}, \eqref{EqGroup:LinearTrilinear},
\eqref{EqGroup:UWBounds}--\eqref{EqGroup:HatUWBounds}, and \eqref{Eq:rMcCormick}--\eqref{EqGroup:UnWnr2UmWmMcCormick} of MBOPF.  Constraint~\eqref{Eq:Rank1Constraint} of RBOPF is relaxed to constraint~\eqref{EqGroup:MinorsLinear} of MBOPF.
	\end{IEEEproof}
Per Propositions~\ref{Proposition:RBOPFRelax} and~\ref{Proposition:MBOPF}, MBOPF~\eqref{EqGroup:MBOPF} is a convex relaxation of the non-convex BOPF~\eqref{EqGroup:OPF}.  The relationship between feasible sets of BOPF, RBOPF, and MBOPF for meshed and radial networks is schematically portrayed in Fig.~\ref{Fig:FeasibleSets}.

\section{Differences With Previous Convex Relaxations}
\label{Sec:LiteratureFormulations}
In this section, we highlight the formulation differences between the proposed approach and previously available convex relaxation techniques for OPF with SVRs.  The premier formulation of~\cite{Robbins2016}, abbreviated here as CIOPF, investigates wye SVRs within the  full SDP relaxation of the  admittance-based  power flow equations combined with following relaxation in place of \eqref{Eq:VoltageGainTrilinear}:
\eq{rCl}{
r_{\min}^2\diag(V_m) \le \diag(V_n) \le r_{\max}^2 \diag(V_m)  \label{Eq:WyeI}
}
Equation~\eqref{Eq:WyeI} can be related to a special case of relaxing \eqref{Eq:Unr2Um} for any $\phi = \phi'$. The work in  \cite{Liu2017}, abbreviated here as CGOPF, similarly uses admittance-based  power flows but employs the chordal SDP relaxation together with
\eq{rCl}{
	r_{\min}^2 V_m \preceq V_n \preceq r_{\max}^2 V_m, \label{Eq:WyeG}
} 
 in place of \eqref{Eq:VoltageGainTrilinear}---based on the simplifying assumption that all SVRs are modeled as gang-operated wye, that is, $r_{nm}^a=r_{nm}^b=r_{nm}^c$ for $(n,m) \in \mc{E}_{\mr{r}}$.  Our previous work~\cite{BazrafshanGatsisZhu2018} also uses the chordal SDP relaxation but includes valid inequalities in the flavor of~\eqref{Eq:WyeI} that are appropriately constructed for closed-delta and open-delta SVRs. However, the applicability of~\cite{BazrafshanGatsisZhu2018} is limited to smaller-sized networks.

 To improve scalability, the formulation MBOPF \eqref{EqGroup:MBOPF} is presented here, where models of  SVRs are incorporated within the \emph{branch flow} formulation of power flow equations. The branch-flow formulation uses the series impedances of transmission lines and transformers, whereas the full or chordal SDP formulations rely on the nodal admittance model. The provided formulation MBOPF is general and  suitable for cases when a mix of wye, closed-delta, and open-delta SVRs are present within the same network. 
 
 If gang operation is desired, the formulation can enforce entries of the effective ratios equal to each other, thereby requiring a single $r_{nm}$ variable for the particular SVR on edge $(n,m) \in \mc{E}_{\mr{r}}$. Such versatility is not available in~\cite{Robbins2016,Liu2017,BazrafshanGatsisZhu2018} as the effective ratio is not an optimization variable in the formulations of the aforementioned works. 
\section{Numerical Experiments}
\label{Sec:NumTests}
The performance of the proposed method is evaluated in this section. Specifically, Section~\ref{Sec:NumTests:Convex} compares the performance of MBOPF to two previously available convex formulations. Section~\ref{Sec:NumTests:Nonlinear} features comparisons with traditional NLP formulations, and provides an instance where MBOPF is preferable over those.

The standard IEEE 13-bus, 37-bus, 123-bus, and 8500-node networks comprising a variety of three-, two-, and one-phase lines are selected for the numerical tests. Transformers are modeled as wye-g--wye-g connections. Switches are replaced by short lines.  Line shunt admittances are ignored, however, capacitors are accounted for as provided by the documentation.

 SVR types for these networks can be  wye, closed-delta, and open-delta. For the 123-bus network, the mixed SVR type means that SVR ID \#1 is modeled as closed-delta, SVR ID \#4 is modeled as open-delta, while the two other SVRs are modeled as wye. For the 8500-node network, the mixed  SVR type means  that SVR IDs \#2 and \#3 are modeled as closed-delta while the two other SVRs are modeled as wye. Voltage regulation on the 8500-node feeder with only open-delta SVRs was not successful, presumably due to lack of
a third tap position, and thus is not reported for any method.

The convex optimization problems are modeled via CVX~\cite{cvx,gb08} and solved by MOSEK~\cite{mosek}.   NLPs are modeled with YALMIP~\cite{Yalmip} and solved by IPOPT~\cite{Wachter2006} through the OPTI interface~\cite{CW12a}.  Experiments in Section~\ref{Sec:NumTests:Convex} are conducted  on a laptop  with a 2-GHz CPU, 8 GB of RAM, and Unix operating system. Experiments in Section~\ref{Sec:NumTests:Nonlinear} are conducted on the same laptop under Microsoft Windows.

\subsection{Performance of the convex relaxation}
\label{Sec:NumTests:Convex}
\begin{table*}[t]
	\scriptsize
	\centering
	\renewcommand{\arraystretch}{1.1}
	\caption{Comparison between various convex formulations}
	\begin{tabular}{|c|c|c|c|c|c|c|c|c|c|c|c|} 
		\hline
		Net.             &  SVR. type & Method               & $c$ (pu)            & $\breve{c}$ (pu)          & $\mr{Gap}$ \%                & $\min \breve{v}$   (pu)              & $\max\breve{v}$  (pu)              & $\breve{v}_{\mr{unb.}}$        &   $ \breve{\Delta}$   (deg.)      & $\lambda_2/\lambda_1$                 & Time (s)           \\     
		\hline \hline 
		\multirow{5}{*}{13-bus} 	&Wye        & CIOPF        & 0.7115      & 0.7135     & 0.2688     & 0.9960     &  1.1118     & 0.0654     &  0.10       & 0.46        &3.21   \\     
		&Wye        & CGOPF         & 0.7140      & 0.7141     & 0.0216     & 0.9911     &  1.0998     & 0.0660     &  0.10       & 0.00        &3.13    \\   
		&Wye        & MBOPF         & 0.7135      & 0.7135     & 0.0033     & 0.9960     &  1.1000     & 0.0534     &  0.10       & 0.00        &3.40        \\
		&Cl.-delta & MBOPF        & 0.7131      & 0.7138     & 0.1038     & 0.9959     &  1.0975     & 0.0435     &  0.86       & 0.00        &3.55     \\  
		&Op.-delta  & MBOPF        & 0.7136      & 0.7145     & 0.1314     & 0.9770     &  1.0967     & 0.0438     &  5.99       & 0.02        &3.41    \\  
		\hline\hline 
		\multirow{5}{*}{37-bus} 	&Wye        & CIOPF         & 1.0351      & 1.0363     & 0.1245     & 0.9236     &  1.0499     & 0.0364     &  1.37       & 0.75        &2.10 \\       
		&Wye        & CGOPF         & 1.0362      & 1.0363     & 0.0147     & 0.9236     &  1.0501     & 0.0363     &  1.37       & 0.00        &2.03     \\   
		&Wye        & MBOPF         & 1.0363      & 1.0363     & 0.0000     & 0.9236     &  1.0501     & 0.0363     &  1.37       & 0.00        & 2.57      \\  
		&Cl.-delta & MBOPF         & 1.0328      & 1.0337     & 0.0912     & 0.9248     &  1.0921     & 0.0294     &  0.35       & 0.00        &3.03   \\     
		&Op.-delta  & MBOPF        & 1.0344      & 1.0347     & 0.0222     & 0.9013     &  1.0938     & 0.0718     &  8.66       & 0.00        &2.74       \\
		\hline \hline
		\multirow{6}{*}{123-bus}		&Wye        & CIOPF       & 0.6997      & 0.7225     & 3.2609     & 0.9592     &  1.0791     & 0.0405     &  2.86       & 0.60        &2.61   \\     
		&Wye        & CGOPF        & 0.7120      & 0.7232     & 1.5806     & 0.9581     &  1.0644     & 0.0424     &  2.89       & 0.00        &2.82     \\
		&Wye        & MBOPF         & 0.7218      & 0.7218     & 0.0026     & 0.9611     &  1.0997     & 0.0307     &  2.89       & 0.00        &3.74\\       
		&Cl.-delta & MBOPF        & 0.7204      & 0.7233     & 0.4061     & 0.9594     &  1.0691     & 0.0491     &  2.13       & 0.01        &5.76 \\     
		&Op.-delta  & MBOPF       & 0.7219      & 0.7224     & 0.0760     & 0.9324     &  1.0989     & 0.0653     &  7.61       & 0.02        &4.33 \\      
		&Mixed      & MBOPF         & 0.7204      & 0.7242     & 0.5315     & 0.9585     &  1.0432     & 0.0449     &  1.08       & 0.02        &4.16       \\
		\hline \hline
		\multirow{5}{*}{8500-node}                	&Wye        & CIOPF         & 0.3972      & 0.4189     & 5.4660     & 0.9407     &  1.1109     & 0.0817     &  5.15       & 0.22        &10.80    \\     
		&Wye        & CGOPF         & 0.4037      & 0.4184     & 3.6370     & 0.9544     &  1.0969     & 0.0693     &  5.03       & 0.00        &13.65    \\       
		&Wye        & MBOPF         & 0.4162      & 0.4176     & 0.3430     & 0.9561     &  1.0999     & 0.0527     &  4.72       & 0.01        &15.23  \\      
		&Cl.-delta & MBOPF         & 0.4157      & 0.4196     & 0.9384     & 0.9106     &  1.0996     & 0.0904     &  5.34       & 0.01        &15.14  \\       
		&Mixed      & MBOPF         & 0.4157      & 0.4195     & 0.9056     & 0.9047     &  1.0961     & 0.0840     &  5.42       & 0.01        & 16.81   \\    
		\hline \hline                                  
	\end{tabular}
	\label{Table:Results}
\end{table*}
 In this section, the OPF cost function is the power import to the distribution network:
\eq{rCl}{c(s_0,\mf{s},V,I)&=& \mr{Re}[\bar{\mb{1}}s_0]. \label{Eq:NumTestsConvexCost}}%
The operational set of power injection [cf.~\eqref{Eq:ConstantPowerSet}] is selected to be a singleton which amounts to the specified load power consumption per phase and per node.
The selection of phase separation parameter $\Delta$ is as follows:    $\Delta=5^{\circ}$ for the wye SVR;   $\Delta=3^{\circ},5^{\circ},10^{\circ},15^{\circ}$ respectively for the 13-bus, 37-bus, 123-bus, and 8500-node networks with closed-delta SVRs;  $\Delta=10^{\circ},10^{\circ}, 15^{\circ}$  for the 13-bus, 37-bus, and 123-bus feeders with open-delta SVRs; and  $\Delta=15^{\circ}$ for the 123-bus and 8500-node networks with mixed SVR types.

After solving the MBOPF, we retrieve the turns ratios using $\breve{r}_{nm}^{\phi}=\sqrt{R_{nm}^{\phi\phi}}$ for wye and closed-delta SVRs.  For open-delta SVRs, we use $\breve{r}_{nm}^{\phi}$ as a solution to the equation $\hat{v}_n-A_{nm}(r_{nm})\hat{v}_m\exp(j \theta_m)=0$ where $\hat{v}_n$ and $\hat{v}_m$ are respectively the spectral decomposition of the rank-1 approximate of $V_n$ and $V_m$ for $(n,m) \in \mc{E}_{\mr{o}}$ and $\theta_m$ is an arbitrary angle variable.   We found this retrieval process for SVR ratios to be more effective in producing feasible voltages during a load-flow, however other methods may also be used.

Upon fixing the ratios, the Z-Bus method is run to obtain voltage solutions $\breve{v}$~\cite{ZBus}. Sufficient conditions for  convergence of the Z-Bus method in three-phase distribution networks are typically satisfied by  IEEE networks~\cite{ZBus,Bernstein2018a}. However, other methods such as the forward-backward sweep may also be used to retrieve voltages~\cite{KerstingBook2001}.  

Table~\ref{Table:Results} provides a summary of performances.   Columns 5--10 respectively provide the following values that are computed based on $\breve{v}$:
\eq{rCl}{ 
\breve{c}&=& \mr{Re}\left[\trace\left(\mr{v}_0 \bar{\breve{v}} \bar{Y}_{0\bullet}\right)\right]\label{Eq:PowerImport}\\
\mr{Gap} \% &=&100 \times (\breve{c}-c)/c \label{Eq:Acc}\\
\min \breve{v}&=& \min_{n, \phi} |\breve{v}_n^\phi| , \quad \max \breve{v} =\max_{n, \phi} |\breve{v}_n^\phi| \label{Eq:MinMaxvLoadFlow}\\
\breve{v}_{\mr{unb.}} &=& \max_{n,\phi}  \left|1- |\breve{v}_n^\phi|/\breve{v}_{\mr{avg}}\right|\label{Eq:vUnbalance} \\
\breve{\Delta} &=& \max\limits_{\phi \in \Omega} |\breve{\theta}_m^{\acute{\phi}}-\breve{\theta}_m^{\phi}-120^\circ|, m: (n,m) \in \mc{E}_{\mr{r}}. \IEEEeqnarraynumspace \label{Eq:MaxAngle}
}

Equation~\eqref{Eq:PowerImport} computes the power import  based on load-flow voltages $\breve{v}$. The notation $Y_{0\bullet}$ denotes the $3 \times (N+1)$ block of the network admittance matrix that corresponds to the slack bus. Equation~\eqref{Eq:Acc} assesses the quality of the objective obtained through the load-flow, that is, $\breve{c}$, in comparison with the objective provided by the corresponding relaxed OPF solution $c$ [cf.~\eqref{Eq:NumTestsConvexCost}]. The quantity $\mr{Gap}$ is the optimality gap, if the load-flow solution $\breve{v}$ turns out to be feasible for the relaxed OPF. The minimum and maximum magnitude of load-flow voltages are given by~\eqref{Eq:MinMaxvLoadFlow}. In~\eqref{Eq:vUnbalance}, $\breve{v}_{\mr{avg}}$ is the average magnitude of voltages.   The quantity $\breve{v}_{\mr{unb.}}$  is a measure of voltage unbalance \cite[eq.~(7.1)]{KerstingBook2001}. Last,  \eqref{Eq:MaxAngle} measures the maximum angle difference from $120^{\circ}$ on the secondary of SVRs based on  $\breve{v}$, assessing validity of Assumption~\ref{Assumption:PhaseSeparation}.

Column 11 of Table~\ref{Table:Results} provides a measure for the rank-1 constraint \eqref{Eq:Rank1Constraint} (and the corresponding rank-1 constraint for the CIOPF and CGOPF formulations) based on the ratio between second-largest eigenvalue ($\lambda_2$) to the largest eigenvalue ($\lambda_1$) of the matrix in \eqref{Eq:Rank1Constraint}  averaged over all non-SVR edges. A  high value of $\lambda_2/\lambda_1$ implies that the matrix is far from being rank-1, while a value close to $0$ implies proximity to a rank-1 solution. Finally, column 12 of Table~\ref{Table:Results} depicts the computation time reported by the solver.
 
We highlight the following key points from Table~\ref{Table:Results}:
 \begin{itemize}
 	\item  For networks with wye SVRs, the optimality gap provided by the proposed MBOPF approach is smaller than the gap obtained from the CIOPF and CGOPF relaxations.   Specifically,  the gap obtained from the proposed approach is below $1\%$ in all networks. The corresponding gap for CIOPF and the CGOPF approaches is above $1\%$ for the 123-bus and 8500-node networks.  
 	
 	\item As a consequence, for networks with only wye SVRs,  the proposed MBOPF approach provides the least-cost feasible solution to the OPF.  Furthermore, the MBOPF yields the smallest voltage unbalance in comparison to CIOPF and CGOPF.
 	
 	\item The proposed MBOPF approach provides a high-quality relaxation for feeders with closed-delta, open-delta, or mixed types of SVRs. In all these cases, the optimality gap is below $1\%$.  In contrast, CIOPF and CGOPF are only valid for networks with wye SVRs. 
 	
 	\item In the IEEE 37-bus feeder, utilizing a closed- or open-delta SVR yields smaller power import costs compared to utilizing a wye SVR, emphasizing the importance of developing convex optimization tools for delta SVRs.
 \end{itemize}
We conclude that  MBOPF  is a reliable and scalable convex formulation for the OPF problem with various types of SVRs.

SVR taps obtained by feeding $\breve{r}_{nm}^{\phi}$ into~\eqref{Eq:Taps} are tabulated in Table~\ref{Table:Taps} for networks with wye SVRs.   Observe that with the exception of the 37-bus feeder,  different formulations of the OPF with  wye SVRs, that is CIOPF, CGOPF, and MBOPF, result in  entirely different tap positions.

\begin{table}[t]
	\centering
	\scriptsize
	\begin{threeparttable}
		\renewcommand{\arraystretch}{1.5}
		\caption{Optimal taps obtained by various convex formulations for wye SVRs}
		\label{Table:Taps}
		\begin{tabular}{cccc}
			SVR ID&  CIOPF &  CGOPF & MBOPF\\
			\hline\hline
			13-1                 & 15,15,15              & 13,13,13             & 15,13,15        \\     
			\hline \hline
			37-1                 & 16,16,16              & 16,16,16             & 16,16,16     \\       
			\hline \hline
			123-1                & 4,-1,-2               & 1,-2,0               & 11,3,7       \\  
			123-2\tnote{*}              & -3                     & 0             &        9      \\  
			123-3\tnote{$\dagger$}            & 1,-4                  & 0,0                  & 11,7     \\               
			123-4                & 15,14,16              & 14,13,13     &         16,16,16            \\          
			\hline \hline
			8500-1               & 7,8,5                 & 8,9,7                & 6,8,6        \\        
			8500-2               & 6,12,10               & 5,4,1                & 5,4,1       \\    
			8500-3               & 3,0,-3                & 4,3,-2               & 12,10,1    \\         
			8500-4               & 1,5,2                 & 1,2,-2               & 6,7,-11      \\  
			\hline    
		\end{tabular}
		\begin{tablenotes}
		\item *  SVR ID \#2 is  single-phase wye
		\item $\dagger$  SVR ID \#3 is a two-phase wye.
	\end{tablenotes}
	\end{threeparttable}
\end{table}

\subsection{Convex vs. NLP formulation}
\label{Sec:NumTests:Nonlinear}

This section compares the solution of the nonlinear BOPF~\eqref{EqGroup:OPF} produced by NLP solvers to that of the convex MBOPF formulation~\eqref{EqGroup:MBOPF}.  In order to highlight the advantages of the convex formulation, we consider an OPF problem  that requires joint optimization of distributed generation (DG) dispatch decisions and  SVR taps.\footnote{The NLP solver for the BOPF formulation managed to find the global optimum to many of the OPF problems of Section~\ref{Sec:NumTests:Convex} upon good initialization. A more complicated  OPF problem is thus presented here to showcase the advantages of a convex formulation over an NLP formulation.}

To this end, the constant-power injection set~\eqref{Eq:ConstantPowerSet} is expressed as $\mc{S}=\prod \mc{S}_m$ where $\mc{S}_m$ is the set of complex constant-power injections $\mf{s}_m$ that satisfy the following constraints for a given load vector $\mf{s}_m^{\mr{l}}$: 
\eq{rCl}{	\label{EqGroup:ConstantPowerSetDg} \IEEEyesnumber \IEEEyessubnumber*
\mf{s}_m &=& \mf{s}_m^{\mr{g}} -\mf{s}_m^{\mr{l}} \label{Eq:sgsl} \\
\left| \mr{Im}\left[\mf{s}_m^{\mr{g}}\right] \right|&\le& \mr{Re}\left[\mf{s}_m^{\mr{g}}\right] \tan\left(\arccos \mr{PF}\right) \label{Eq:PF} \\
\IEEEeqnarraymulticol{3}{c}{\sqrt{\mr{Re}\left[\mf{s}_m^{\mr{g}}\right]^2+\mr{Im}\left[\mf{s}_m^{\mr{g}}\right]^2} \le \mf{s}_{\max}} \label{Eq:Capacity} 
	}
In~\eqref{EqGroup:ConstantPowerSetDg}, $\mf{s}_m^{\mr{g}}$ is a variable representing the complex power generation of the DG at node $m$, and the constants $\mf{s}_{\max}$ and $\mr{PF}$ respectively  denote the apparent power capacity and the maximum power factor (capacitive of inductive) of the DG.
The objective is to minimize the total amount of real power injection and power import to the distribution network, i.e., 
\eq{rCl}{
c(s_0,\mf{s},v,i)=c(s_0,\mf{s}^{\mr{g}})=\mr{Re}\left[\bar{\mb{1}}s_0+\bar{\mb{1}}\mf{s}^{\mr{g}}\right], \label{Eq:NumTestsObjNLP}
}
where $\mf{s}^{\mr{g}}=\{\mf{s}_m^{\mr{g}}\}_{m \in \mc{N}}$ collects all load vectors. 

Voltage limits are set to $v_{\min}=0.95$ and $v_{\max}=1.05$. DG is only connected to buses with three available phases and specific values of $\mr{PF}=0.9$ and $\mf{s}^{\max}=0.001$ pu  have been selected. The BOPF is initialized with voltages obtained from a load-flow when the network contains no DGs and the SVR taps are set to $0$.  BOPF~\eqref{EqGroup:OPF} and  MBOPF~\eqref{EqGroup:MBOPF} are then solved by the respective solvers and system configurations as detailed at the beginning of Section~\ref{Sec:NumTests}. 
 
The resulting tap ratios $\breve{r}$ and constant-power injections $\mf{s}$  from the optimization stage are fed into a Z-Bus load-flow to compute the feasible objective
\eq{rCl}{
\breve{c}&=& \mr{Re}\left[\trace\left(\mr{v}_0 \bar{\breve{v}} \bar{Y}_{0\bullet}\right)+\bar{\mb{1}}\mf{s}^{\mr{g}}\right].\label{Eq:FeasibleObjWithDG}
}
The gap between the optimal value of the optimization stage and the feasible objective $\breve{c}$, that is the discrepancy between~\eqref{Eq:NumTestsObjNLP} and~\eqref{Eq:FeasibleObjWithDG}, as well the minimum and maximum voltages are computed similar to~\eqref{Eq:Acc} and~\eqref{Eq:MinMaxvLoadFlow}.

Table~\ref{Table:BNLPResultsSg} tabulates the performance of the convex MBOPF and NLP BOPF formulations for the  8500-node network with mixed SVRs.  Table~\ref{Table:Taps2} reveals  that the taps obtained from MBOPF and BOPF can be significantly different.  A first glance on Table~\ref{Table:BNLPResultsSg} shows that the NLP formulation BOPF remarkably finds a solution whose cost is practically equal to the cost obtained by the convex MBOPF. However, column 6 of Table~\ref{Table:BNLPResultsSg} reveals that solving the convex formulation MBOPF is significantly faster---a remarkable speedup of at least 20 times.  

It is worth emphasizing that the good performance of the NLP formulation BOPF is in general dependent on the  initialization point. Furthermore, NLP solvers potentially could get stuck in a local minimum. On the contrary, the convex MBOPF formulation provides a useful lower bound along with a feasible solution within a very short time-span. Accordingly, the solution obtained by MBOPF can be powerful for both assessing the quality of solution of NLP solvers and for warm-starting them.

\begin{table}[t]
	\scriptsize
	\centering
	\renewcommand{\arraystretch}{1.1}
	\caption{Convex vs. nonlinear programming}
	\begin{tabular}{|c|c|c|c|c|c|} 
		\hline 
	   Method                        & $\breve{c}$ (pu)          & $\mr{Gap}$ \%                & $\min \breve{v}$   (pu)              & $\max\breve{v}$  (pu)            & Time (s)           \\    
	   \hline \hline
	MBOPF & 0.3871 & 0.0020 & 0.9609 & 1.0415 & 8.71 \\
	BOPF & 0.3872 & 0.0000& 0.9733 & 1.0112 & 187.17 \\	            
		\hline 
	\end{tabular}
	\label{Table:BNLPResultsSg}
\end{table}

\begin{table}[t]
	\centering
	\scriptsize
	\begin{threeparttable}
		\renewcommand{\arraystretch}{1.2}
		\caption{Optimal taps obtained by MBOPF and BOPF}
		\label{Table:Taps2}
		\begin{tabular}{ccc}
			SVR ID&  MBOPF & BOPF\\
			\hline\hline
			8500-1                             & -1,-1,6              & -3,-3,1       \\        
			8500-2\tnote{*}                           & 3,1,1            & 1,1,1      \\    
			8500-3\tnote{*}                           & -1,1,-4             & -1,1,4   \\         
			8500-4                             & -3,-2,-10              & 0,2,-4    \\  
			\hline    
		\end{tabular}
		\begin{tablenotes}
			\item *  Closed-delta SVRs
		\end{tablenotes}
	\end{threeparttable}
\end{table}

\section{Concluding Remarks}
\label{Sec:Conclusions}

This paper introduces an SDP framework for OPF problems that include tap selection of the most common SVRs in practice. The branch flow model of the power flow equations is adopted and extended to handle SVR edges. A phase separation assumption is introduced, which is realistic and adopted only for the secondary voltages of SVRs.  Specialized techniques are developed to relax the various non-convexities that show up due to the rigorous modeling of SVRs. The resultant convex program represents a quite tight relaxation that is coupled with a tap recovery scheme leading to very small optimality gaps even in large-scale networks. Future work includes extending the present framework to multi-period OPF and tap selection problems that limit the cycling and wear-and-tear of the regulation equipment.

\color{black}
\appendices
\section{Useful linear algebra results} \label{Appendix:UsefulLinearAlgebra}
\begin{lemma} \label{Lemma:DiagToHermitian}
	For two complex vectors $u$ and $w$, we have that 
	\eq{rCl}{\diag(u) w^*=\diag(u\bar{w}). \label{Eq:SimpleIdentity}}
\end{lemma} 
\begin{IEEEproof}
	The proof is omitted due to its simplicity. 
\end{IEEEproof}
\section{Proof of Proposition~\ref{Proposition:Retrieve}} \label{Appendix:RetrieveProof}
	The proof follows the  procedure in \cite{Gan2014} but extends it to handle SVR edges. The proof is based on induction. At the $n$-th iteration with $\mc{N}_{\mr{visit}}^{(n)}$, $v_n$ is given and the following holds:
	\eq{rCl}{V_n=v_n \bar{v}_n, n \in \mc{N}_{\mr{visit}}^{(n)}. \label{Eq:Vnvnvnbar}}
	We have to then show that \eqref{Eq:Ohms}, \eqref{EqGroup:SVRGains}, \eqref{Eq:NetPowerInjectionBreakDownNonConvex}, \eqref{Eq:PowerBalanceNonConvex}, and \eqref{Eq:VoltageBounds} are satisfied. First notice that~\eqref{Eq:Ohms} and \eqref{Eq:VoltageGain} are satisfied by construction of Algorithm~\ref{Algorithm:Retrieve}.  To prove that \eqref{Eq:NetPowerInjectionBreakDownNonConvex} is satisfied, we will show that \eqref{Eq:AuxV} holds for $(n,m) \in \mc{E}$ and $m \notin \mc{N}_{\mr{visit}}^{(n)}$.  The equality~\eqref{Eq:AuxV} and constraint \eqref{Eq:VBounds} automatically yield \eqref{Eq:VoltageBounds}. To prove that \eqref{Eq:PowerBalanceNonConvex} is satisfied, we will show that \eqref{Eq:AuxI} and \eqref{Eq:AuxPower} hold for $(n,m) \in \mc{E}$ and the diagonal of \eqref{Eq:AuxPrime} holds for $(n,m) \in \mc{E}_{\mr{r}}$, that is \eqref{Eq:DiagofAuxPower} holds.
	Finally, based on \eqref{Eq:VoltageGain},  \eqref{Eq:AuxPower}, \eqref{Eq:DiagofAuxPower}, and \eqref{Eq:SVRPowers} Lemma~\ref{Lemma:PowerSubstituteCurrent} proves that \eqref{Eq:CurrentGain} also holds.   Therefore, it suffices to show that from $v_n$ satisfying \eqref{Eq:Vnvnvnbar}, we can construct $(v_m, i_{nm}, i'_{nm})$ that satisfy \eqref{Eq:AuxV}--\eqref{Eq:AuxPower} and \eqref{Eq:DiagofAuxPower}.  
	
	For every $(n,m) \in \mc{E}$, \eqref{Eq:PDConstraint} and \eqref{Eq:Rank1Constraint} hold, which implies
	\eq{rCl}{ \bmat{V_n & S_{nm} \\ \bar{S}_{nm} & I_{nm}}&=& \bmat{u \\ w}\bmat{\bar{u} & \bar{w}},\label{Eq:SpectralDecompose}}
	for some  complex vectors $u$ and $w$. Therefore, 
	\eq{rCl}{\label{EqGroup:SpectralDecomposePerEntry} \IEEEyesnumber \IEEEyessubnumber*
		V_n &=& u\bar{u} \label{Eq:Vnuubar} \\
		S_{nm}&=& u \bar{w} \label{Eq:Snmuwbar} \\
		I_{nm}&=& w\bar{w} \label{Eq:Inmwwbar} }
	Equations~\eqref{Eq:Vnuubar} together with~\eqref{Eq:Vnvnvnbar} imply that 
	\eq{rCl}{
		v_{n}&=& u\exp(j \theta) \label{Eq:vnutheta}
	}
	for some vector $\theta$ and the product in \eqref{Eq:vnutheta} is entrywise.  Using~\eqref{Eq:Vnuubar},~\eqref{Eq:Snmuwbar} and~\eqref{Eq:vnutheta} in Algorithm~\ref{Algorithm:Retrieve} update \ref{Eq:inmUpdate} yields
	\eq{rCl}{
		i_{nm} &=& \frac{1}{\bar{u} u} w\bar{u} u \exp(j \theta)= w\exp(j\theta). \label{Eq:inmwtheta}
	}
	Substituting $u$ and $w$ in \eqref{Eq:Snmuwbar} and \eqref{Eq:Inmwwbar} readily yield \eqref{Eq:AuxI} and \eqref{Eq:AuxPower}.    To obtain \eqref{Eq:AuxV}, if $(n,m) \in \mc{E}_{\mr{t}}$ then
	\eq{rCl}{v_{m}\bar{v}_m =(v_n-Z_{nm}i_{nm})(\bar{v}_n - \bar{i}_{nm} \bar{Z}_{nm}) \hfill \notag \\
		= V_n+Z_{nm}I_{nm}\bar{Z}_{nm}  -(S_{nm}\bar{Z}_{nm}+Z_{nm}\bar{S}_{nm})=V_m \IEEEeqnarraynumspace \label{Eq:TrlinevmvmbarVm}} where the last equality comes from  \eqref{Eq:OhmsMatrix}.
	If $(n,m) \in \mc{E}_{\mr{r}}$, then Algorithm~\ref{Algorithm:Retrieve} update \ref{Eq:vmSVRUpdate} gives
	\eq{rCl}{v_{m}\bar{v}_m &=& A_{nm}^{-1} v_n \bar{v}_n \bar{A}_{nm}^{-1} = A_{nm}^{-1} V_n \bar{A}_{nm}^{-1}=V_m \label{Eq:SVRvmvmbarVm}}
	where the last equality comes from  \eqref{Eq:VoltageGainTrilinear}.  Therefore,  \eqref{Eq:AuxV}--\eqref{Eq:AuxPower} hold. It remains to show that \eqref{Eq:DiagofAuxPower} holds for $(n,m) \in \mc{E}_{\mr{r}}$.   From Algorithm~\ref{Algorithm:Retrieve} update \ref{Eq:inmPrimeSVRUpdate} it holds that 
	\eq{rCl}{
		\diag(S'_{nm})&=&\diag(v_m) (i'_{nm})^* =\diag(v_m \bar{i}_{nm}) \label{Eq:DiagSnmPrimevminmPrimeBar}
	}
	where the last equality uses Lemma \ref{Lemma:DiagToHermitian}.  Lemma~\ref{Lemma:PowerSubstituteCurrent} can now be invoked to show that \eqref{Eq:CurrentGain} also holds. \hfill \IEEEQED
 
 \section{Proof of Proposition~\ref{Proposition:UWBounds}} \label{Appendix:ProofUWBounds}
 		To prove \eqref{EqGroup:UWBounds}, notice from
 		\eqref{Eq:AuxV} that $V_m^{\phi\phi'}=v_m^{\phi} \bar{v}_m^{\phi'}$. Therefore, for diagonal elements it holds that $U_m^{\phi\phi}=\Re{V_m^{\phi\phi}}=|v_m^{\phi}|^2$ and $W_m^{\phi \phi}=\Im{V_m^{\phi\phi}}=0$ which implies that we have $U_{\min}^{\phi \phi}=v_{\min}^2$, $U_{\max}^{\phi\phi}=v_{\max}^2$, while $W_{\min}^{\phi \phi}=W_{\max}^{\phi\phi}=0$. For the $(\phi, \acute{\phi})$-th element it holds that $U_m^{\phi\acute{\phi}}=\Re{V_m^{\phi\acute{\phi}}}=|v_m^{\phi}| |v_m^{\acute{\phi}}| \cos(\phi- \acute{\phi})$ and $W_m^{\phi \acute{\phi}}=|v_m^{\phi}| |v_m^{\acute{\phi}}| \sin(\phi-\acute{\phi})$ which together with \eqref{Eq:ThetaBounds} implies that we have $U_{\min}^{\phi \acute{\phi}}=v_{\max}^2 \cos\left(120^\circ+\Delta\right)$, 	$U_{\max}^{\phi \acute{\phi}}=v_{\min}^2 \cos\left(120^\circ-\Delta\right)$, and $W_{\min}^{\phi \acute{\phi}}=v_{\min}^2 \sin\left(120^\circ+\Delta\right)$ and $W_{\max}=v_{\min}^2 \sin\left(120^\circ-\Delta\right)$. The remaining entries are filled by acknowledging that $V_m$ is Hermitian. 
 		
 		Bounds in \eqref{EqGroup:TildeUWBounds} and \eqref{EqGroup:HatUWBounds} are computed next. For wye SVRs, it holds that $D_{nm}=I$ and $F_{nm}=\mb{O}$.  Hence, for wye SVRs,  $\tilde{U}_{\min}$, $\tilde{U}_{\max}$, $\tilde{W}_{\min}$, and $\tilde{W}_{\max}$ are respectively equal  to $U_{\min}$, $U_{\max}$, $W_{\min}$,  and $W_{\max}$ while $\hat{U}_{\min}$, $\hat{U}_{\max}$, $\hat{W}_{\min}$, and $\hat{W}_{\max}$ are  zeros. 
 		For closed-delta and open-delta SVRs, the expressions for $\tilde{U}_{\min}$, $\tilde{U}_{\max}$, $\tilde{W}_{\min}$,  $\tilde{W}_{\max}$, $\hat{U}_{\min}$, $\hat{U}_{\max}$, $\hat{W}_{\min}$, and $\hat{W}_{\max}$ contain more terms. However, it turns out that the $(\psi, \psi')$-th element of $\tilde{U}_{nm}$ and $\hat{U}_{nm}$, denoted here by $u^{\psi\psi'}$,  is of the form 
 			\eq{rCl}{
 		u^{\psi \psi'}&=& \sum\limits_{\phi \in \Omega}\left[ a_{\phi \acute{\phi}}^{\psi \psi'} |v_m^{\phi}| |v_m^{\acute{\phi}}|  \cos\left(\phi - \acute{\phi}\right) \right. \notag \\
 		&& \left. - \:  b_{\phi \acute{\phi}}^{\psi \psi'} |v_m^{\phi}| |v_m^{\acute{\phi}}| \cos\left(\phi - \acute{\phi}\right) \right] \notag \\
 		&&+ \: \sum\limits_{\phi \in \Omega}c_{\phi}^{\psi \psi'}|v_m^{\phi}|^2 - d_{\phi}^{\psi \psi'}|v_m^{\phi}| ^2 \IEEEeqnarraynumspace \label{Eq:upsipsi'}}
 
\noindent where $a_{\phi\acute{\phi}}^{\psi \psi'}$, $b_{\phi\acute{\phi}}^{\psi \psi'} $, $c_{\phi}^{\psi \psi'}$, and $d_{\phi \phi}^{\psi \psi'}$ are all non-negative constants for $\phi, \psi, \psi' \in \Omega$. 
Therefore, based on \eqref{Eq:VoltageBounds} and \eqref{Eq:ThetaBounds}, the bounds on $u^{\psi\psi'}$ are  given by
	\eq{rCl}{\IEEEyesnumber \label{EqGroup:upsipsi'Bounds} \IEEEyessubnumber*
	u_{\min}^{\psi \psi'}&=& \sum\limits_{\phi \in \Omega}\left[ a_{\phi \acute{\phi}}^{\psi \psi'} v_{\max}^2 \cos\left(120^\circ+\Delta\right) \right. \notag \\
	&& \left. - \:  b_{\phi \acute{\phi}}^{\psi \psi'} v_{\min}^2 \cos\left(120^\circ-\Delta\right) \right] \notag \\
	&& + \: \sum\limits_{\phi \in \Omega}c_{\phi}^{\psi \psi'}v_{\min}^2 - d_{\phi}^{\psi \psi'}v_{\max}^2 \label{Eq:uminpsipsi'} \\
	u_{\max}^{\psi \psi'}&=& \sum\limits_{\phi \in \Omega}\left[ a_{\phi \acute{\phi}}^{\psi \psi'} v_{\min}^2 \cos\left(120^\circ-\Delta\right) \right. \notag \\
&& \left. - \:  b_{\phi \acute{\phi}}^{\psi \psi'} v_{\max}^2 \cos\left(120^\circ+\Delta\right) \right] \notag \\
&&  + \: \sum\limits_{\phi \in \Omega}c_{\phi}^{\psi \psi'}v_{\max}^2 - d_{\phi}^{\psi \psi'}v_{\min}^2 \label{Eq:umaxpsipsi'} \IEEEeqnarraynumspace}
Similarly,  it turns out that the $(\psi, \psi')$-th element of $\tilde{W}_{nm}$ and $\hat{W}_{nm}$, denoted here by $w^{\psi\psi'}$,  is of the form
\eq{rCl}{
	u^{\psi \psi'}&=& \hspace{-1mm}\sum\limits_{\phi \in \Omega}\left[ e_{\phi \acute{\phi}}^{\psi \psi'} |v_m^{\phi}| |v_m^{\acute{\phi}}|  \sin\left(\phi - \acute{\phi}\right)  \notag \right. \\
	 && \left. -\: f_{\phi \acute{\phi}}^{\psi \psi'} |v_m^{\phi}| |v_m^{\acute{\phi}}| \sin\left(\phi - \acute{\phi}\right) \right] 
	 \label{Eq:wpsipsi'} \IEEEeqnarraynumspace }where $e_{\phi\acute{\phi}}^{\psi \psi'}$ and $f_{\phi\acute{\phi}}^{\psi \psi'} $ are non-negative constants. Hence, based on \eqref{Eq:VoltageBounds} and \eqref{Eq:ThetaBounds}, the bounds on $w^{\psi\psi'}$ are  given by
	\eq{rCl}{\IEEEyesnumber \label{EqGroup:wpsipsi'Bounds} \IEEEyessubnumber*
	w_{\min}^{\psi \psi'}&=& \hspace{-1mm} \sum\limits_{\phi \in \Omega}\left[ e_{\phi \acute{\phi}}^{\psi \psi'} v_{\min}^2 \sin\left(120^\circ+\Delta\right) \right. \notag \\ 
	&& \left. - \:  b_{\phi \acute{\phi}}^{\psi \psi'} v_{\max}^2 \sin\left(120^\circ-\Delta\right) \right] \label{Eq:wminpsipsi'} \\
	w_{\max}^{\psi \psi'}&=& \hspace{-1mm}\sum\limits_{\phi \in \Omega}\left[ e_{\phi \acute{\phi}}^{\psi \psi'} v_{\max}^2 \sin\left(120^\circ-\Delta\right) \right.  \notag \\
	&& \left. -\:  f_{\phi \acute{\phi}}^{\psi \psi'} v_{\min}^2 \sin\left(120^\circ+\Delta\right) \right]. \IEEEeqnarraynumspace \label{Eq:wmaxpsipsi'} }

\bibliographystyle{IEEEtran}
\IEEEtriggeratref{18}
\bibliography{vr}

\end{document}